\documentclass[final,1p,number,sort&compress]{siamart1116}

\usepackage[english]{babel}
\usepackage{graphicx,epstopdf,epsfig}
\usepackage{amsfonts,epsfig,fancyhdr,graphics, hyperref,amsmath,amssymb}

\numberwithin{equation}{section}
\usepackage{mathtools,mathdots}
\usepackage{algorithm}
\usepackage{algpseudocode}
\usepackage{amssymb}
\usepackage{amsmath}
\usepackage{booktabs,url}

\newcommand{\Names}{Christian Bertram, Heike Fa\ss bender}
\newcommand{\Title}{A rational Krylov subspace based projection method for solving CAREs}

\newtheorem{remark}[theorem]{Remark}

\renewcommand{\theequation}{\thesection.\arabic{equation}}
\renewtheorem{theorem}{Theorem}[section]

\numberwithin{equation}{section} 

\begin{document}

\bibliographystyle{plain}

\setcounter{page}{1}

\thispagestyle{empty}

 \title{A class of Petrov-Galerkin Krylov methods for algebraic Riccati equations}

\author{
Christian Bertram\thanks{Institute for Numerical Analysis, TU Braunschweig,
Braunschweig, Germany.}
\and
Heike Fa\ss bender\thanks{Institute for Numerical Analysis, TU Braunschweig,
Braunschweig, Germany (h.fassbender@tu-braunschweig.de).}
}

\pagestyle{myheadings}
\markboth{\Names}{\Title}

\maketitle

\begin{abstract}
A class of (block) rational Krylov subspace based projection method for solving large-scale continuous-time algebraic Riccati equation (CARE) 
$0 = \mathcal{R}(X) := A^HX + XA + C^HC - XBB^HX$ with a large, sparse $A$ and $B$ and $C$ of full low rank is proposed.
The CARE is projected onto a block rational Krylov subspace $\mathcal{K}_j$ spanned by blocks of the form $(A^H - s_kI)^{-1}C^H$ for some shifts $s_k, k = 1, \ldots, j.$  The considered projections do not need to be orthogonal and are built from the matrices appearing in the block rational Arnoldi decomposition associated to $\mathcal{K}_j.$ The resulting projected Riccati equation is solved for the small square Hermitian $Y_j.$ Then the Hermitian low-rank approximation $X_j = Z_jY_jZ_j^H$ to $X$  is set up where the columns of $Z_j$  span $\mathcal{K}_j.$ 
The residual norm $\|R(X_j )\|_F$ can be computed efficiently via the norm of a readily available $2p \times 2p$ matrix. We suggest to reduce the rank of the approximate solution $X_j$ even further by truncating small eigenvalues from $X_j.$ This truncated approximate solution can be interpreted as the solution of the Riccati residual projected to a subspace of $\mathcal{K}_j.$ This gives us a way to efficiently evaluate the norm of the resulting residual.  Numerical examples are presented. 
\end{abstract}




\section{Introduction}
The numerical solution of large-scale algebraic Riccati equations
\begin{equation}\label{ric0}
0 = \mathcal{R}(X) := A^HX+XA+C^HC-XBB^HX 
\end{equation}
with a large, sparse matrix $A \in \mathbb{C}^{n \times n},$ and matrices $B\in \mathbb{C}^{n \times m}$ and $C\in \mathbb{C}^{p \times n}$ is of interest in a number of applications as noted in \cite{BenS13,LinS15,SimSM14} and references therein. As usual, we are interested in finding the stabilizing solution $X;$ that is, the  Hermitian positive semidefinite solution $X =X^H$ such that the spectrum of the closed loop matrix $(A-BB^HX)$ lies in the open left half plane $\mathbb{C}_-.$
This solution exists and is unique \cite{BinIM11,LanR95} when $(A,B)$ is stabilizable (that is, rank$[A-zI,B]=n$ for each value of $z$ in the closed right half plane) and $(A,C)$ is detectable (that is, 
$(A^H,C^H)$ is stabilizable). These conditions are generically fulfilled and are assumed to hold in the remainder.
For our discussion, it is  further assumed that $B$ and $C$ have full column and row rank, resp.,  with $m,p \ll n.$ Then
 it is well-known that the sought-after solution $X$ will often have a low numerical rank (that is, its numerical rank $q$ is $\ll n$) \cite{BenB16}. Thus, it can be expressed in the form $X = ZYZ^H$ for some full rank matrices $Z \in \mathbb{C}^{n \times q}, Y=Y^T \in \mathbb{C}^{q \times q}$ with $q < n.$ This allows for the construction of iterative methods that approximate $X$ with a series of low-rank matrices stored in low-rank factored form. There are several methods (e.g., rational Krylov subspace methods, low-rank Newton-Kleinman methods and Newton-ADI-type methods)
which produce such a low-rank approximation; see, e.g. \cite{AmoB10,BenBKS18,BenHSW16,HeyJ09,LinS15,Pal19,Sim16,SimSM14,WonB05,WonB07} and \cite{BenBKS20} for an overview.). 

Here we are concerned with projection type methods which project the original Riccati equation \eqref{ric0} onto an appropriate subspace, impose a Galerkin condition on the projected problem and expand its solution back to the whole space.
This idea can be summarized as follows (see, e.g., \cite{Sim16,BenBKS20}). Assume a sequence of nested subspaces $\mathcal{N}_k \subseteq \mathcal{N}_{k+1}, k \geq 1$ is generated with $\dim\mathcal{N}_k =kp.$ Let $Q_k\in \mathbb{C}^{n \times kp}$ denote the matrix whose columns vectors are an orthonormal basis of $\mathcal{N}_k.$ Then the Galerkin condition reads
\[
Q_k^H \mathcal{R}(X_k) Q_k = 0.
\]
This gives a small-scale Riccati equation for $Y_k=Y_k^H \in \mathbb{C}^{kp \times kp}$
\begin{equation}\label{smallVk}
A_k^HY_k+Y_kA_k-Y_kB_kB_k^HY_k+C_k^HC_k = 0
\end{equation}
with $A_k = Q_k^HAQ_k\in \mathbb{C}^{kp \times kp}, B_k = Q_k^HB\in \mathbb{C}^{kp \times m}$ and $C_k^H=Q_k^HC^H\in \mathbb{C}^{kp \times p}.$ 

If $(A_k,B_k)$ and $(A_k^H,C_k^H)$ are stabilizable and detectable, resp., then the existence of a stabilizing solution of \eqref{smallVk} is ensured. In \cite[Proposition 3.3]{Sim16}, a general condition on $A$ and $B$ is given which
ensures stabilizability of $(A_k,B_k)$. Unfortunately, in practice it is difficult to check whether this condition holds.
But as the set of matrices $(\tilde{A},\tilde{B},\tilde{C})$ with $(\tilde{A},\tilde{B})$ stabilizable
and $(\tilde{A},\tilde{C})$ detectable
is dense in $\mathbb{C}^{\tilde{n }\times \tilde{n}}\times \mathbb{C}^{\tilde{n} \times m} \times \mathbb{C}^{p\times \tilde{n}}$, most likely \eqref{smallVk} will have a unique stabilizing Hermitian positive semidefinite solution.
Thus, in general, a unique Hermitian stabilizing solution $Y_k$ of this small-scale Riccati equation exists, such that the approximate solution
\[
X_k = Q_kY_kQ_k^H
\] 
can be constructed. When $k = n$, $\mathcal{R}(X_k) =0$ must hold and the exact solution has been found. The effectiveness of this approach depends on the choice of $\mathcal{N}_k$. The approximation spaces explored in the literature are all based on block (rational) Krylov subspaces, see \cite{DruS11,HeyJ09,LinS15,Pal19,Sim16,SimSM14}. For a short historic review of these ideas (stemming usually from algorithms to solve Lyapunov equations) and a list of references  see \cite[Section 2]{SimSM14} or \cite[Section 2]{BenBKS20}.
In \cite{HeyJ09} the matrix $Q_k$ is constructed as a basis of the extended block Krylov subspace $\kappa_k(A^H,C^H) + \kappa_k(A^{-H},A^{-H}C^H)$ \cite{DruK98,KniS10} for the standard block Krylov subspace 
\begin{equation}\label{eq_standardKrylov}
\kappa_k(A^H,C^H) = \operatorname{range}([C^H, A^HC^H, (A^H)^2C^H, \ldots, (A^H)^{k-1}C^H]).
\end{equation}
In \cite{DruS11,SimSM14} in addition $Q_k$ is also generated from the block rational Krylov subspace  \cite{Ruh84}
\begin{equation}\label{eq_rationalKrylov}
\mathfrak{K}_{k}\coloneqq \kappa_{k}(A^H,C^H,\mathcal{S}_{k-1})= \operatorname{range}([C^H, (A^H-s_1I)^{-1}C^H, \ldots, (A^H-s_{k-1}I)^{-1}C^H])
\end{equation}
 for a set of shifts $\mathcal{S}_{k-1} = \{s_1, \ldots, s_{k-1}\} \subset \mathbb{C}$ disjoint from the spectrum $\Lambda(A^H).$ The use of the rational Krylov subspace $\mathfrak{K}_{k}$ is explored further in \cite{Sim16}. The resulting projection algorithm is usually termed RKSM algorithm. In \cite{BenBKS20}, it was observed that  projecting to rational Krylov subspaces often fared better than projecting to extended Krylov subspaces. 


In both \cite{LinS15} and \cite{BenBKS18} methods are proposed which can be interpreted as projecting the large-scale Riccati equation \eqref{ric0} onto the block rational  Krylov subspace $\mathcal{K}_k$
\begin{equation}\label{eq_ratKryohneC}
\mathcal{K}_k \coloneqq \mathcal{K}_k(A^H,C^H,\mathcal{S}_k) =\operatorname{range}([(A^H-s_1I)^{-1}C^H, \ldots, (A^H-s_kI)^{-1}C^H])
\end{equation}
 for a set of shifts $\mathcal{S}_k = \{s_1, \ldots, s_k\} \subset \mathbb{C}$ disjoint from the spectrum $\Lambda(A^H).$
In \cite{LinS15} approximate solutions to the CARE \eqref{ric0} are constructed by running subspace iterations
on the Cayley transforms of the Hamiltonian matrix $\left[ \begin{smallmatrix} A & BB^H\\ C^HC& -A^H\end{smallmatrix}\right]\in \mathbb{C}^{2n \times 2n}$ associated to \eqref{ric0}. It is observed that the columns of the orthonormal matrix $Q_k$ in the resulting approximate solution $X_k^\text{cay} = Q_kY_k^\text{cay}Q_k^H$ span $\mathcal{K}_k$. Moreover, it is proven  in \cite[Theorem 4.4]{LinS15} that the subspace iteration yields (under certain conditions) a matrix $Q_k$ such that besides \eqref{eq_ratKryohneC} also $Q_k^H\mathcal{R}(X_k^\text{cay})Q_k =0$ holds. In case the reduced problem \eqref{smallVk} admits a unique stabilizing solution, this gives equivalence between the subspace iteration and the Galerkin projection onto the rational Krylov subspace $\mathcal{K}_k.$ The RADI algorithm suggested in \cite{BenBKS18}  is a generalization of the Cholesky-factored variant of the Lyapunov ADI method. It generates approximations $X_k^\text{radi} =  Z_kY_k^\text{radi}Z_k^H$ where the columns of $Z_k$ can be interpreted as a nonorthogonal basis of $\mathcal{K}_k$. As noted in \cite[Section 4.2]{BerF23},  $X_j^\text{radi}$ can be interpreted as the solution of a projection of the large-scale Riccati equation \eqref{ric0} onto the Krylov subspace $\mathcal{K}_k$ employing an oblique projection. 
%
For both methods, the rank of the Riccati residual $\mathcal{R}(X_k)$ is always equal to $p.$ The norm of $\mathcal{R}(X_k)$ can be evaluated efficiently giving a cheap stopping criterion.

We take up the idea of employing a Petrov-Galerkin condition for the Riccati residual equation \eqref{ric0}.
In our discussion we solely concentrate on projecting onto a block rational Krylov subspace $\mathcal{K}_k$ \eqref{eq_ratKryohneC} which so far has not been considered in the literature as an idea on its own. We make use of the (generalized) block rational Arnoldi decomposition (BRAD) $A^HV_{k+1}\underline{K}_k = V_{k+1} \underline{H}_k$ \cite[Definition 2.2]{ElsG20} where $\operatorname{range}(V_{k+1})=  \mathfrak{K}_{k+1}.$ 
As $\operatorname{range}(V_{k+1}\underline{K}_k)=  \mathcal{K}_{k}$, choosing $Z_k = V_{k+1}\underline{K}_k$ and a suitable matrix $W_k$ yields a (not necessarily orthogonal) projection $\Pi_k = Z_k(W_k^HZ_k)^{-1}W_k^H$ onto  $\mathcal{K}_{k}.$ 
The projected CARE $\Pi_k\mathcal{R}(X_k)\Pi_k^H=0$ gives a small-scale Riccati equation whose coefficient matrices are given essentially by the matrices of the BRAD.
 The resulting projected Riccati equation is solved for the small square Hermitian $Y_k.$ 
Then the Hermitian low-rank approximation $X_k = Z_kY_kZ_k^H$ to $X$  is set up. 
Details of this approach are explained in the Section \ref{sec3}. As explained in Section \ref{sec4}, 
the residual norm $\|\mathcal{R}(X_k )\|_F$ can be computed efficiently via the norm of a readily available $2p \times 2p$ matrix avoiding the explicit computation of the large-scale matrices $X_k $ and $\mathcal{R}(X_k ).$ A secondary result of that discussion is that rank$(\mathcal{R}(X_k ))=2p,$ unlike for the RADI algorithm where the residual is of rank $p.$
In Section \ref{sec5}, we suggest to reduce the rank of the approximate solution $X_k$ even further by truncating small singular values from $X_k.$ This truncated approximate solution can be interpreted as the solution of the Riccati residual projected to a subspace of $\mathcal{K}_k.$ As a result, we obtain a way to efficiently evaluate the norm of the resulting residual. The transfer of the results to the generalized Riccati  equation is explained in Section \ref{sec6}.  Numerical examples are presented in Section \ref{sec7}. The results discussed here can be found in slightly different form in the PhD thesis \cite{Ber21}.


\section{Preliminaries}\label{sec:p1}
Rational Krylov spaces were initially proposed by Ruhe in the 1980s for the purpose of solving large sparse eigenvalue
problems \cite{Ruh84,Ruh94,Ruh98}, more recent work on block Krylov subspaces by Elsworth and G\"uttel can be found in \cite{ElsG20}. We briefly recall some definitions and results on block rational Krylov subspaces and block rational Arnoldi decompositions from \cite[Section 1 and 2]{ElsG20}.

Given  a set $\mathcal{S}_k =\{s_1, \ldots, s_k \} \subset \mathbb{C}\backslash\{\infty\}$ disjoint from the spectrum $\Lambda(A^H).$   In case $A,B,C$ in \eqref{ric0} are real matrices, any complex $s_k$ should to be accompanied by its complex-conjugate partner $\overline{s}_j,$ so that the set $\mathcal{S}_k$ is closed under conjugation. In that case, all algorithms suggested here can be implemented in real arithmetic.

In the following, we will consider the block Krylov subspaces \eqref{eq_rationalKrylov} and \eqref{eq_ratKryohneC}, that is, 
\begin{align*}
\mathfrak{K}_{k+1} &= \kappa_{k+1}(A^H,C^H,\mathcal{S}_{k})= \operatorname{range}([C^H, (A^H-s_1I)^{-1}C^H, \ldots, (A^H-s_{k}I)^{-1}C^H]),\\\mathcal{K}_k &= \mathcal{K}_k(A^H,C^H,\mathcal{S}_k) =\operatorname{range}([(A^H-s_1I)^{-1}C^H, \ldots, (A^H-s_kI)^{-1}C^H]).
\end{align*}
 The columns of the starting block $C^H$ lie in $\mathfrak{K}_{k+1},$ but not in $\mathcal{K}_k$.
 It is assumed that the $(k+1)p$ columns of $\mathfrak{K}_{k+1}$ are linearly independent. 

There is a one-to-one correspondence between block rational Krylov spaces and so-called block rational Arnoldi decompositions.
In \cite{Sim16} the decomposition
\begin{equation}\label{ratarn_sim}
A^HQ_{j+1} = Q_{j+2}\underline{T}_{j+1}^H
\end{equation}
is used where the columns of $Q_{j+1} \in \mathbb{C}^{n \times (j+1)p}$ are orthonormal ($Q_{j+1}^HQ_{j+1} = I_{(j+1)p}$) and span either the extended block Krylov subspace or the block rational Krylov subspace $\mathfrak{K}_{j+1}.$ We will consider the more general relation of the form
\begin{equation}\label{rad}
A^HV_{j+1}\underline{K}_j = V_{j+1} \underline{H}_j
\end{equation}
 called a (generalized) orthonormal block rational Arnoldi decomposition (BRAD)  as introduced in \cite{ElsG20}.
The following conditions hold
\begin{itemize}
\item  the columns of $V_{j+1} \in \mathbb{C}^{n \times (j+1)p}$ are orthonormal such that $\operatorname{range}(V_{j+1})=  \mathfrak{K}_{j+1}$,
\item the matrices $\underline{H}_j$ and $\underline{K}_j$ are block upper Hessenberg matrices of size $(j+1)p\times jp$ where at least one of the matrices  $H_{i+1,i}$ and $K_{i+1,i}$ is nonsingular,
\end{itemize}
\begin{equation}\label{eq_HK}
 \underline{H}_j =\left[ \begin{smallmatrix}
H_{11} & H_{12} & \cdots & H_{1,j-1} & H_{1j}\\
H_{21} & H_{22} & \cdots & H_{2,j-1} & H_{2j}\\
& \ddots & \ddots & \vdots& \vdots\\
& & H_{j-1,j-2} & H_{j-1,j-1} & H_{j-1,j}\\
& & & H_{j,j-1} & H_{jj}\\
&&&&H_{j+1,j}
\end{smallmatrix} \right], 
\underline{K}_j = \left[\begin{smallmatrix}
K_{11} & K_{12} & \cdots & K_{1,j-1} & K_{1j}\\
K_{21} & K_{22} & \cdots & K_{2,j-1} & K_{2j}\\
& \ddots & \ddots & \vdots& \vdots\\
& & K_{j-1,p-2} & K_{j-1,j-1} & K_{j-1,j}\\
& & & K_{j,j-1} & K_{jj}\\
&&&&K_{j+1,j}
\end{smallmatrix}\right],
\end{equation}
 \begin{itemize}
\item $\beta_iK_{i+1,i} =\gamma_iH_{i+1,i}$ with scalars $\beta_i,\gamma_i \in \mathbb{C}$ such that $|\beta_i|+|\gamma_i|\neq 0$ for all $i = 1, \ldots, j,$ and
\item  the quotients $\beta_i/\gamma_i$, $i = 1, \ldots, j,$ called poles of the BRAD, correspond to the shifts $s_1, \ldots, s_j \in \mathcal{S}_j.$
\end{itemize}
An algorithm which constructs an orthonormal BRAD can be found in \cite[Algorithm 2.1]{ElsG20}, see also \cite{BerEG14}.

Due to \cite[Lemma 3.2 (iii)]{ElsG20}, both matrices $\underline{K}_j$ and  $\underline{H}_j$ are of full rank $jp$. This does not imply that all subdiagonal blocks $H_{i+1,i}$ and $K_{i+1,i}$ are nonsingular.
As noted in \cite[Lemma 3.2 (ii)]{ElsG20}, in case one of the subdiagonal blocks $H_{i+1,i}$ and $K_{i+1,i}$ is singular, it is the zero matrix. Thus, in case $H_{i+1,i}=0$ in \eqref{rad}, the matrix $K_{i+1,i}$ must be nonsingular. As $\beta_iK_{i+1,i} =\gamma_iH_{i+1,i}$ has to hold with $|\beta_i|+|\gamma_i|\neq 0$, this implies $\beta_i=0.$  Hence, this can only happen for a zero shift in $\mathcal{S}_j.$ In case $K_{i+1,i}=0$ in \eqref{rad}, the matrix $H_{i+1,i}$ must be nonsingular. As $\beta_iK_{i+1,i} =\gamma_iH_{i+1,i}$ has to hold with $|\beta_i|+|\gamma_i|\neq 0$, this implies $\gamma_i=0$ and a shift at infinity. Due to our assumption that the roots are finite, this case is excluded here. Thus, all blocks $K_{i+1,i}$ are nonsingular, while some blocks $H_{i+1,i}$ may be singular.

Any BRAD can be transformed into an equivalent one of the form 
\[
A^H\check{V}_{j+1}\underline{\check{K}}_j = \check{V}_{j+1} \underline{\check{H}}_j
\]
with $\underline{\check{K}}_j = \left[ \begin{smallmatrix} 0_{p\times jp}\\I_{jp}\end{smallmatrix}\right],$ $\operatorname{range}(\check{V}_{j+1}) = \mathfrak{K}_{j+1},$ $\check V_{j+1} \left[ \begin{smallmatrix} I_p\\0_{jp \times p}\end{smallmatrix}\right]=  C^H\Phi, \Phi \in \mathbb{C}^{p\times p}$, and an upper Hessenberg matrix $\underline{\check {H}}_j$, see \cite[Section 2]{BerF23}. 
This implies
\begin{align}
\operatorname{range}( V_{j+1}\underline{{K}}_j) = \operatorname{range}(\check V_{j+1}\underline{\check{K}}_j) &= \mathcal{K}_j, \label{range_info1}\\
\operatorname{range}( V_{j+1}\underline{{H}}_j) = \operatorname{range}(\check V_{j+1}\underline{\check{H}}_j)&= A^H\mathcal{K}_j. \label{range_info2}
\end{align}

\begin{remark}\label{rem1}
Consider the BRAD $A^HV_{j+1}\underline{K}_j = V_{j+1} \underline{H}_j.$ 
Let the thin QR decomposition $\underline{K}_j = QR$ be given with an orthonormal upper Hessenberg matrix $Q\in \mathbb{C}^{(j+1)p \times jp}$ and an upper triangular $R\in \mathbb{C}^{jp \times jp}.$ As $\underline{K}_j$ is of full rank, $R$ is nonsingular. With $\hat{\underline{K}}_j = Q$ and $\hat{\underline{H}}_j = \underline{H}_jR^{-1}$ we obtain the equivalent BRAD
\[
A^HV_{j+1}\hat{\underline{K}}_j = V_{j+1} \hat{\underline{H}}_j
\]
with an orthonormal $\hat{\underline{K}}_j.$ In case $V_{j+1}$ is orthonormal as well, the matrix $Z_j = V_{j+1}\hat{\underline{K}}_j $ is orthonormal.
\end{remark}

\section{General Projection Method}\label{sec3}
In this section we present in detail our approach to reduce the Riccati residual \eqref{ric0} to a low dimensional one using a Petrov-Galerkin projection.
In order to do so, let the columns of $Z_j  \in \mathbb{C}^{n \times jp}$ span  some approximation space $\mathcal{K}_j$.
We choose a matrix $W_j \in \mathbb{C}^{n\times jp}$ such that its columns span the test space $\mathcal{L}_j$ and $W_j^HZ_j \in \mathbb{C}^{jp \times jp}$ is nonsingular. Thus,
\[
\Pi_j \coloneqq Z_j(W_j^HZ_j)^{-1} W_j^H \in \mathbb{C}^{n \times n}
\]
is a  (not necessarily orthogonal) projection onto $\mathcal{K}_j$ along $\mathcal{L}_j^\perp.$ 
Assume  that the solution $X$ of \eqref{ric0} can be approximated by
\begin{equation}\label{eqX}
X_j = Z_jY_jZ_j^H
\end{equation}
 for some Hermitian matrix $Y_j \in \mathbb{C}^{jp \times jp}.$ Then, imposing a Petrov-Galerkin condition on $\mathcal{R}(X_j) $
\begin{align*}
0 &= \Pi_j\mathcal{R}(X_j) \Pi_j^H\\
&= Z_j\left\{ \left((W_j^HZ_j)^{-1}W_j^HA^H Z_j\right ) Y_j + Y_j \left(Z_j^HAW_j(W_j^HZ_j)^{-H}\right)
 \right.\\
&~~~~~~ \left.+ \left( (W_j^HZ_j)^{-1} W_j^H C^H\right) \left(C W_j(W_j^HZ_j)^{-H}\right) + Y_j (Z_j^HB)  (B^HZ_j)Y_j\right\} Z_j^H
\end{align*}
one obtains a small scale Riccati equation for $Y_j$
\begin{equation}\label{redric}
0 = A_j^HY_j+Y_jA_j +C_j^HC_j-Y_jB_j B_j^HY_j
\end{equation}
with $A_j = Z_j^HAW_j(W_j^HZ_j)^{-H} \in \mathbb{C}^{jp\times jp},$  $C_j = C W_j(W_j^HZ_j)^{-H}\in \mathbb{C}^{p \times jp},$ and $ B_j = Z_j^HB \in \mathbb{C}^{jp \times m}.$

In particular, we suggest to project onto the block rational Krylov subspace $ \mathcal{K}_j $ as in \eqref{eq_rationalKrylov}; that is, $\mathcal{K}_j=\mathcal{K}_j(A^H,C^H,\mathcal{S}_j).$ A suitable basis for this space is constructed via the generalized orthonormal RAD
\begin{equation}\label{rado}
A^HV_{j+1}\underline{K}_j = V_{j+1} \underline{H}_j,
\end{equation}
where the columns of $V_{j+1}\underline{K}_j$ span the space $ \mathcal{K}_j=\mathcal{K}_j(A^H,C^H,\mathcal{S}_j)$  (see \eqref{range_info1}). Thus, we choose
\begin{equation}\label{eq_Z}
Z_j \coloneqq V_{j+1}\underline{K}_j  \in \mathbb{C}^{n \times jp}.
\end{equation}
Note that the columns of $Z_j$ are in general not orthonormal.
As $\underline{K}_j$ is of full rank, it is possible to choose a matrix $\underline{L}_j \in \mathbb{C}^{(j+1)p\times jp}$ such that $\underline{L}_j^H\underline{K}_j \in \mathbb{C}^{jp \times jp}$ is nonsingular. Set $W_j \coloneqq V_{j+1}\underline{L}_j \in \mathbb{C}^{n \times jp}.$
Denote (as before) the space spanned by the columns of the matrix $W_j$ by $\mathcal{L}_j.$ Hence,
$\Pi_j  Z_j(W_j^HZ_j)^{-1} W_j^H \in \mathbb{C}^{n \times n}$
is  a projection onto $\mathcal{K}_j$ along $\mathcal{L}_j^\perp.$
Noting that $W_j^HZ_j = \underline{L}_j^HV_{j+1}^HV_{j+1}\underline{K}_j = \underline{L}_j^H\underline{K}_j$ holds and making use of \eqref{rado}, the coefficient matrices in \eqref{redric} can be expressed as
\begin{align}\label{redABC}
\begin{split}
A_j &\coloneqq \underline{H}_j^H \underline{L}_j(\underline{K}_j^H\underline{L}_j)^{-1}    \in \mathbb{C}^{jp \times jp},\\
B_j &\coloneqq \underline{K}_j^HV_{j+1}^HB \in \mathbb{C}^{jp \times m},\\
C_j &\coloneqq CV_{j+1}\underline{L}_j(\underline{K}_j^H\underline{L}_j)^{-1} \in \mathbb{C}^{p\times jp}.
\end{split}
\end{align}
The solution $Y_j$ of \eqref{redric} fully determines the approximate solution $X_j = V_{j+1}\underline{K}_j Y_j \underline{K}_j^HV_{j+1}$ \eqref{eqX} of \eqref{ric0}.

The resulting algorithm is summarized in Algorithm \ref{alg1}.  In the first step $V_1$ is chosen from the thin QR decomposition of $C^H = V_1R.$  This yields $CV_{j+1} = \left[ \begin{smallmatrix}R^H & 0_{p\times jp}\end{smallmatrix}\right].$ Once the iteration is stopped, the approximate solution $X_j = V_{j+1}\underline{K}_j Y_j \underline{K}_j^HV_{j+1}^H$ can be set up (theoretically). In a large scale setting, one may consider only its low-rank factor $V_{j+1}\underline{K}_jG_j$ with a Cholesky(-like) decomposition $Y_j =G_jG_j^H.$

\begin{algorithm}
\caption{General projection method for solving \eqref{ric0}}\label{alg1}
\begin{algorithmic}[1]
\Require System matrices $A, B, C$  and set of shifts $\mathcal{S}_j = \{s_1, \ldots, s_j\} \subset \mathbb{C}\backslash\{\infty\} $ with $\mathcal{S}_j \cap \Lambda(A^H)= \emptyset.$ 
\Ensure Approximate solution $X_j= V_{j+1}\underline{K}_j Y_j \underline{K}_j^HV_{j+1}^H.$
\State Compute thin QR decomposition $C^H = V_1 R$ with $R\in \mathbb{C}^{p \times p}.$
\State Set $j = 1.$
\While {not converged}
\State  Obtain next shift $\mu$ from $\mathcal{S}_j.$
\State Expand orthonormal BRAD to obtain  $V_{j+1} = \left[V_j ~~\hat V_{j+1}\right]$ orthonormal, $\underline{K}_j$ and $\underline{H}_j.$
\State Choose $\underline{L}_j.$
\State Compute $A_j = \underline{H}_j^H \underline{L}_j(\underline{K}_j^H\underline{L}_j)^{-1}.$
\State Compute $B_j = \underline{K}_j^HV_{j+1}^HB$.
\State Compute $C_j = \left[R^H~~0_{p \times jp}\right] \underline{L}_j(\underline{K}_j^H\underline{L}_j)^{-1}.$
\State Solve $A_j^HY_j + Y_j A_j +C_j^HC_j - Y_jB_jB_j^HY_j =0.$
\State Compute $\left\| \mathcal{R}(X_j)\right\|$ (see Algorithm \ref{alg2}).
\State Set $j = j+1.$
\EndWhile
\end{algorithmic}
\end{algorithm}

From the rank conditions on the matrices involved in \eqref{redABC} it follows that $A_j$ is nonsingular,  and the matrices $C_j$ and $B_j$ are of full rank as $\operatorname{rank}C_j =p$ and $\operatorname{rank}B_j = \min\{pj,m\}.$ Whether \eqref{redric} has a unique stabilizing solution can not be read off in general. But as the set of matrices $(A_j,B_j,C_j)$ with $(A_j,B_j)$ stabilizable and $(A_j,C_j)$ detectable is dense in $\mathbb{C}^{jp \times jp}\times \mathbb{C}^{jp \times m} \times \mathbb{C}^{p\times jp}$, most likely \eqref{redric} will have a unique stabilizing Hermitian positive semidefinite solution. Only in very rare cases of our numerical experiments the small-scale equation \eqref{redric} could not be solved for a stabilizing solution.

The effectiveness of the proposed projection method depends on the choice of the set $\mathcal{S}_j$ of shifts and of the space $\mathcal{L}_j$/the matrix $\underline{L}_j.$ 
Here we briefly mention three  different choices of $\underline{L}_j:$ 
\begin{itemize}
\item A natural choice is $\underline{L}_j \coloneqq \underline{K}_j$, which yields a Galerkin projection as $\mathcal{L}_j = \mathcal{K}_j.$
\item Another choice is $\underline{L}_j \coloneqq \underline{H}_j$, which yields a Petrov-Galerkin projection as $\mathcal{L}_j = A^H\mathcal{K}_j$ due to \eqref{range_info2}.
\item A more general choice is  $\underline{L}_j \coloneqq \alpha \underline{H}_j - \beta \underline{K}_j$ with $\alpha, \beta \in \mathbb{C}$ and $|\alpha|+|\beta| \neq 0$ which yields a Petrov-Galerkin projection in case $\alpha \neq 0.$ 
\end{itemize}

\begin{remark}
Assume that an orthonormal BRAD with orthonormal $\underline{K}_j$ (see Remark \ref{rem1}) is used.
Then the matrix $Q_j \coloneqq V_{j+1}\underline{K}_j$ has orthonormal columns. Its columns are an orthonormal basis of the Krylov subspace $\mathcal{K}_j.$
Set $\underline{L}_j  = \underline{K}_j,$ so the projection is orthogonal. Then it follows from \eqref{redABC}
\begin{eqnarray*}
A_j &=& \underline{H}_j^H \underline{K}_j (\underline{K}_j^H \underline{K}_j)^{-1} 
  =  \underline{H}_j^H V_{j+1}^HV_{j+1}\underline{K}_j  =  \underline{K}_j^H V_{j+1}^HA V_{j+1}\underline{K}_j = Q_j^HAQ_j,\\
B_j &=& \underline{K}_j^HV_{j+1}^HB = Q_j^HB,\\
C_j &=& C Q_j
\end{eqnarray*}
due to $\underline{K}_j^H \underline{K}_j=I_j,$ $V_{j+1}^HV_{j+1}=I_{j+1}$ and the BRAD \eqref{rado}. These reduced matrices look like the ones in \eqref{smallVk}. Recall that here we project onto the Krylov subspace $\mathcal{K}_j$, while in most of the literature \cite{HeyJ09,LinS15,Pal19,Sim16,SimSM14} a projection on the Krylov subspace $\mathfrak{K}_j$ is used.
In numerical experiments no advantage concerning accuracy when using $Q_j^HAQ_j, Q_j^HB,CQ_j$ has been observed.
\end{remark}
\begin{remark} Assume that $V_{j+1}$ is some matrix with orthonormal columns.
If $Z_j = V_{j+1}\underline{K}_j$ is chosen with some arbitrary $\underline{K}_j$ not satisfying \eqref{rado}, then $A_j$ in \eqref{redABC} is given by  $A_j = \underline{K}_j^H V_{j+1}^HAV_{j+1}\underline{L}_j(\underline{K}_j^H\underline{L}_j)^{-1}.$  The derivations in the Sections \ref{subsec_eff}, \ref{subsec_shifts},  \ref{sec4} and \ref{sec5} do not hold in this case, as they depend on \eqref{rado}. 
\end{remark}

{\subsection{Efficient Implementation}\label{subsec_eff}
Algorithm \ref{alg1}  should not be implemented as formulated above. The block upper Hessenberg form of the matrices $\underline{K}_j$ and $\underline{H}_j$ should be utilized so that the matrices $A_j, B_j, C_j$ are not explicitly set up in each step, but rather obtained by a (simple) update from $A_{j-1}, B_{j-1}, C_{j-1},$ resp.. Moreover, in step 7 and 9 of Algorithm \ref{alg1}, a linear system of equations with $jp$, resp. $p$, right hand sides has to be solved for the same $jp \times jp$ system matrix. This would require roughly $(jp)^3$ flops if the structure of the system matrix cannot be utilized.

Let us first consider the computation of $B_j.$ As $\underline{K}_j$ is a block upper Hessenberg matrix which grows by one block row/column in each iteration step we have from \eqref{eq_HK}
\begin{equation}\label{eq_K}
\underline{K}_j = \begin{bmatrix} \underline{K}_{j-1} & k_{j}\\ 0_{p\times (j-1)p} & K_{j+1,j}\end{bmatrix} \text{\quad with \quad}
k_{j}=\begin{bmatrix} K_{1j}^H ~~ K_{2j}^H~~\cdots~~K_{jj}^H\end{bmatrix}^H \in \mathbb{C}^{jp\times p}.
\end{equation}
Thus, as $V_{j+1}= \begin{bmatrix}V_j~~\hat V_{j+1}\end{bmatrix}$, it holds that
\begin{equation}\label{eq_B}
B_j = \underline{K}_j^HV_{j+1}^HB =  \left[ \begin{smallmatrix} \underline{K}_{j-1}^H & ~0_{(j-1)p\times p}\\ k_{j}^H & K_{j+1,j}^H\end{smallmatrix} \right]
\left[ \begin{smallmatrix} V_{j}^H\\ \hat V_{j+1}^H\end{smallmatrix}\right] B = \left[\begin{smallmatrix}B_{j-1} \\k_{j}^HV_j^HB+ K_{j+1,j}^H\hat V_{j+1}^HB\end{smallmatrix}\right].
\end{equation}
Thus, in each iteration step, $B_j$ can be updated cheaply from $B_{j-1}=\underline{K}_{j-1}^HV_{j}^HB $  as only the last block row of $B_j$ needs to be determined, the rest of $B_j$ is already known. Moreover, the re-computation of $V_j^HB$ for the last block row can be avoided by storing $\tilde B_j = V_{j+1}^HB = [ \tilde B_{j-1} ~~\hat V_{j+1}^HB]$ in each iteration step. Then $B_j = [k_{j}^H ~~ K_{j+1,j}^H]\tilde B_j.$

Next, we consider the computation of the matrix products $\underline{H}_j^H\underline{L}_j$ and $\underline{K}_j^H\underline{L}_j.$ If in each iteration step  a new $\underline{L}_j$ is selected without taking the previous $\underline{L}_{j-1}$ into account, these products have to be computed from scratch in each iteration step. This is different in case $\underline{L}_j$ is chosen as a block upper Hessenberg matrix whose leading part is given by $\underline{L}_{j-1}$ like $\underline{K}_j$ or $\underline{H}_j$ in \eqref{eq_HK} (see also \eqref{eq_K}),
\begin{equation}\label{eq_L}
\underline{L}_j= \begin{bmatrix} \underline{L}_{j-1} & \ell_{j}\\ 0_{p \times (j-1)p} & L_{j+1,j}\end{bmatrix}.
\end{equation}
Then  only the last block row and column of $\underline{H}_j^H\underline{L}_j$ and $\underline{K}_j^H\underline{L}_j$ need to be computed in case  $\underline{H}_{j-1}^H\underline{L}_{j-1}$ and $\underline{K}_{j-1}^H\underline{L}_{j-1}$ are already known, e.g.,
\[
\underline{K}_j^H\underline{L}_j= \left[ \begin{smallmatrix} \underline{K}_{j-1}^H & ~&0_{(j-1)p\times p}\\ k_{j}^H & &K_{j+1,j}^H\end{smallmatrix} \right]
\left[ \begin{smallmatrix} \underline{L}_{j-1} &~& \ell_{j}\\ 0_{p \times (j-1)p} & ~&L_{j+1,j}\end{smallmatrix} \right] =
\left[ \begin{smallmatrix} \underline{K}_{j-1}^H \underline{L}_{j-1} & ~& \underline{K}_{j-1}^H\ell_{j}\\ k_{j}^H \underline{L}_{j-1} &~&k_{j}^H\ell_{j}+  K_{j+1,j}^HL_{j+1,j}\end{smallmatrix} \right].
\]

We need to be able to solve linear systems with the coefficient matrix $\underline{K}_j^H\underline{L}_j$ efficiently in order to set up $A_j$ and $C_j$.
One idea is to make use of the block LDU decomposition
\begin{equation}\label{eq_ldu}
\underline{K}_j^H\underline{L}_j= 
\begin{bmatrix} I_{(j-1)p} & 0_{(j-1)p\times p}\\ \Omega_j& I_p\end{bmatrix}
\begin{bmatrix} \Psi_j & 0_{(j-1)p\times p}\\ 0_{p \times (j-1)p} & \Gamma_j\end{bmatrix}
\begin{bmatrix} I_{(j-1)p} &  \Delta_j\\0_{p \times (j-1)p}  & I_p\end{bmatrix}
\end{equation}
with
\begin{align*}
\Omega_j &= k_j^H\underline{L}_{j-1}(K_{j-1}^H\underline{L}_{j-1})^{-1}\in \mathbb{C}^{p \times (j-1)p}, \\
\Delta_j &= (\underline{K}_{j-1}^H \underline{L}_{j-1})^{-1} \underline{K}_{j-1}^H\ell_j\in \mathbb{C}^{(j-1)p \times p},\\
\Psi_j &= \underline{K}_{j-1}^H\underline{L}_{j-1}\in \mathbb{C}^{(j-1)p \times (j-1)p},\\
\Gamma_j &= k_{j}^H\ell_{j}+  K_{j+1,j}^HL_{j+1,j} - 
k_j^H\underline{L}_{j-1}(\underline{K}_{j-1}^H \underline{L}_{j-1})^{-1} \underline{K}_{j-1}^H\ell_j \\
&=k_{j}^H\ell_{j}+  K_{j+1,j}^HL_{j+1,j} -  k_j^H\underline{L}_{j-1}\Delta_j \in \mathbb{C}^{p \times p}.
\end{align*}
Setting up this decomposition, two linear systems with the system matrix $\Psi_j=\underline{K}_{j-1}^H\underline{L}_{j-1}$ need to be solved in order to determine $\Omega_j$ and $\Delta_j$, both with $p$ right hand sides. For $\underline{K}_{j-1}^H\underline{L}_{j-1}$, there is a block LDU decomposition just like \eqref{eq_ldu} which may be used to solve those systems recursively. In particular, for $\Delta_j=\left[\begin{smallmatrix}\Delta_{j,1}\\\Delta_{j,2}\end{smallmatrix}\right], \Omega_j = \left[ \Omega_{j,1}~~\Omega_{j,2}\right]$ and $\underline{K}_{j-1}^H \ell_j= \left[\begin{smallmatrix}\mathfrak{h}_{1}\\\mathfrak{h}_{2}\end{smallmatrix}\right], k_j^H\underline{L}_{j-1} = \left[\mathfrak{g}_1~~\mathfrak{g}_2\right]$ with $\Delta_{j,1}, \mathfrak{h}_1, \mathfrak{g}_1^T \in \mathbb{C}^{(j-2)p \times p}$ and
$\Delta_{j,2}, \mathfrak{h}_2, \mathfrak{g}_2^T \in \mathbb{C}^{p\times p}$ 
we obtain
\begin{align*}
\Delta_{j,2} &= \Gamma_{j-1}^{-1}(\mathfrak{h}_2 - \Omega_{j-1}\mathfrak{h}_1 ),\\
\Delta_{j,1} &= \Psi_{j-1}^{-1}\mathfrak{h}_1 - \Delta_{j-1}\Delta_{j,2},\\
\Omega_{j,2} &= (\mathfrak{g}_2 -\mathfrak{g}_1\Delta_{j-1})\Gamma_{j-1}^{-1}\\
\Omega_{j,1} &= \mathfrak{g}_1\Psi_{j-1}^{-1} -\Omega_{j,2}\Omega_{j-1}.
\end{align*}
Thus, continuing in this fashion, in order to determine $\Delta_j$ and $\Omega_j$ $2j$ linear systems with the $p \times p$ coefficient matrices $\Gamma_\ell, \ell = 1, \ldots, j$ need to be solved for $p$ right hand sides each.  The computational cost for each of these steps is dominated by computing the LU decomposition of  $\Gamma_\ell$ which costs $\mathcal{O}(p^3)$ flops. Hence, setting up the LDU decomposition of $\underline{K}_j^H\underline{L}_j$ amounts to $\mathcal{O}(jp^3)$ flops. 

Now we will concentrate on computing $A_j=\underline{H}_j^H\underline{L}_j(\underline{K}_j^H\underline{L}_j)^{-1}$ and assume that we have already determined
\[
 \widetilde{A}_j=\underline{H}_j^H\underline{L}_j =
 \begin{bmatrix} 
\underline{H}_{j-1}^H \underline{L}_{j-1} &  \underline{H}_{j-1}^H\ell_{j}\\ 
h_{j}^H \underline{L}_{j-1} &h_{j}^H\ell_{j}+  H_{j+1,j}^HL_{j+1,j}
\end{bmatrix} =\begin{bmatrix} (\widetilde A_j)_{11} & (\widetilde A_j)_{12}\\ (\widetilde A_j)_{21} & (\widetilde A_j)_{22}\end{bmatrix}.
\]
The matrix 
\[
A_j =\underline{H}_j^H\underline{L}_j(\underline{K}_j^H\underline{L}_j)^{-1}= \widetilde{A}_j (\underline{K}_j^H\underline{L}_j)^{-1} = \begin{bmatrix} (A_j)_{11} & (A_j)_{12}\\ (A_j)_{21} & (A_j)_{22}\end{bmatrix} 
\]
is given by
\begin{align*}
 (A_j)_{12}&=\left( (\widetilde A_j)_{12} - (\widetilde A_j)_{11}\Delta_j\right)\Gamma_j^{-1}, &
(A_j)_{11} =(\widetilde A_j)_{11}\Psi_j^{-1}-(A_j)_{12}\Omega_j, \\
 (A_j)_{22} &= \left((\widetilde A_j)_{22} - (\widetilde A_j)_{21})\Delta_j\right)\Gamma_j^{-1},&
(A_j)_{21} = (\widetilde A_j)_{21}\Psi_j^{-1}- (A_j)_{22}\Omega_j,
\end{align*}
due to \eqref{eq_ldu} where $(\widetilde A_j)_{11}\Psi_j^{-1} = \underline{H}_{j-1}^H \underline{L}_{j-1} (\underline{K}_{j-1}^H \underline{L}_{j-1} )^{-1}$ is already known from the previous iteration. Thus, in order to determine $A_j$ two linear systems with the $p \times p$ coefficient matrix $\Gamma_j$ need to be solved for $(j-1)p$ right hand sides each.  Computing the LU decomposition of $\Gamma_j$ costs $\mathcal{O}(p^3)$ flops, the $(2j-2)p$ forward and backward solves amount to $\mathcal{O}(jp^3)$ flops. 
Computing $(\widetilde A_j)_{21}\Psi_j^{-1}$  can be done in the same way as the calculation of $\Omega_j$ discussed above. The LU decompositions of the $\Gamma_\ell. \ell = 1, \ldots, j$ should be reused here such that only additional forward and backward solve are needed here.

Finally, we consider the computation of $C_j = C V_{j+1}\underline{L}_j(\underline{K}_j^H\underline{L}_j)^{-1}.$ 
With $R^H=CV_1$, we obtain $CV_{j+1} = \left[ \begin{smallmatrix}R^H & 0_{p\times jp}\end{smallmatrix}\right].$
Assuming that $\underline{L}_j$ is chosen as a block upper Hessenberg matrix  as in \eqref{eq_L}, it follows that 
\begin{align*}
\tilde C_j = CV_{j+1}\underline{L}_j &= \begin{bmatrix}R^H & 0_{p\times jp}\end{bmatrix} 
\left[\begin{smallmatrix}
L_{11} & L_{12} & \cdots & L_{1,j-1} & L_{1j}\\
L_{21} & L_{22} & \cdots & L_{2,j-1} & L_{2j}\\
& \ddots & \ddots & \vdots& \vdots\\
& & L_{j-1,j-2} & L_{j-1,j-1} & L_{j-1,j}\\
& & & L_{j,j-1} & L_{jj}\\
&&&&L_{j+1,j}
\end{smallmatrix}\right]\\
 &=
\begin{bmatrix}
R^HL_{11} & R^HL_{12} & \cdots & R^HL_{1,j-1} & R^HL_{1j}
\end{bmatrix}\\
&= \begin{bmatrix}CV_{j}\underline{L}_{j-1}&  R^HL_{1j} \end{bmatrix}  = \begin{bmatrix}\tilde C_{j-1} &  (\tilde C_j)_{12} \end{bmatrix}.
\end{align*}
Thus, $CV_{j+1}\underline{L}_j$ can be updated cheaply from $CV_{j}\underline{L}_{j-1}$ as only the block  $R^HL_{1j}$ needs to be determined, the rest of $CV_{j+1}\underline{L}_j$ is already known. We obtain $C_j =   \tilde C_j  (\underline{K}_j^H\underline{L}_j)^{-1} = [ (C_j)_{11}~~(C_j)_{12}]$ with $(C_j)_{12}\in \mathbb{C}^{p \times p}$ and \eqref{eq_ldu} via
\begin{align*}
 (C_j)_{12}&=  ((\tilde C_j)_{12} - \tilde C_{j-1}\Delta_j)\Gamma_j^{-1}, & (C_j)_{11} &= \tilde C_{j-1}(\underline{K}_{j-1}^H\underline{L}_{j-1})^{-1}-(C_j)_{12}\Omega_j.
\end{align*}
The matrix $ \tilde C_{j-1}(\underline{K}_{j-1}^H\underline{L}_{j-1})^{-1}$ is already known from the previous step. Thus, in order to determine $C_j$ from $C_{j-1},$  there is just  one linear system with the $p \times p$ system matrix $\Gamma_j$ and $p$ right hand sides to be solved. Thus, the cost for setting up $C_j$ is  dominated by the cost for the LU decomposition of $\Gamma_j,$ that is, $\mathcal{O}(p^3)$ flops.

In summary, each iteration step can be implemented such that just  $\mathcal{O}(jp^3)$ flops are needed in case $\underline{L}_j$ is chosen as a block upper Hessenberg matrix whose leading part is given by $\underline{L}_{j-1}$ like $\underline{K}_j$ or $\underline{H}_j$ in \eqref{eq_HK} (see also \eqref{eq_K}). This approach implies that all matrices $\Omega_\ell, \Delta_\ell^T \in \mathbb{C}^{p \times (\ell-1)p}$ and $\Gamma_\ell \in \mathbb{C}^{p \times p}, \ell = 1, \ldots, j$ have to be stored.  Otherwise, that is, without taking the previous $\underline{L}_{j-1}$ into account, each iteration step will cost  $\mathcal{O}((jp)^3)$ flops.


\subsection{Shift Selection}\label{subsec_shifts}
For fast convergence the choice of the poles of the block rational Krylov subspace used in the approximate solution $X_j$ is crucial. 
It would be desirable to determine a set of shifts $\mathcal{S}_j$ such that the approximation $X_j^\ast$ satisfies $\|X-X_j^\ast\| = \min_{\mathcal{S}_j}\|X-X_j\|$ in some norm where $X_j$ has been computed by Algorithm \ref{alg1} using $\mathcal{S}_j.$ 
This question has been considered for symmetric Sylvester equations \cite{BenB14}. Related work on Lyapunov equations can be found in \cite{Wol14,VanV10}. Upto now, it is an open question on how to compute optimal shifts in the sense described above for Riccati equations.

Many shift strategies based on other ideas exist, but their description is beyond the scope of this work. In Section 5 of \cite{DruS11,Sim16} the shift selection for methods projecting \eqref{ric0} onto the rational Krylov subspace $\mathfrak{K}_j$ \eqref{eq_rationalKrylov} is discussed.
The subspace iteration method proposed in \cite{LinS15} can (under certain assumptions) be interpreted as projecting \eqref{ric0} onto the rational Krylov subspace $\mathcal{K}_j$ \eqref{eq_ratKryohneC}. Hence, the comments on the shift selection given in \cite{LinS15} hold for the approach proposed in this paper. 
Moreover, the subspace iteration method is equivalent to the RADI algorithm \cite{BenBKS18} (when using $X_0=0$ and the same set of shifts).
A detailed discussion of the shift selection in the context of the RADI algorithm and related, equivalent  methods is given in \cite[Section 4.5]{BenBKS18}. In particular, the idea of choosing $s_{j+1}$ in order to minimize $\|\mathcal{R}(X_{j+1})\|$ once $X_j$ is fixed is pursued \cite[Section 4.5.2]{BenBKS18}. A recent numerical comparison of different solvers of \eqref{ric0} including the choice of shifts can be found in 
 \cite{BenBKS20}.

\section{Efficient Residual Norm Evaluation}\label{sec4}
Typically, the residual norm $\left\|\mathcal{R}(X_j)\right\|$ is used as an indicator whether $X_j$ is a good approximation to the desired solution $X$ of \eqref{ric0}. In the previous section we replaced solving the large-scale $n \times n$ problem \eqref{ric0} by solving of a much smaller one \eqref{redric}. In order to be able to use  $\left\|\mathcal{R}(X_j)\right\|$ as a convergence indicator, an equivalent small-scale expression has to be found as 
$\mathcal{R}(X_j) \in \mathbb{C}^{n \times n}$ can not  be computed explicitly due to storage constraints.
In the following we will derive such an expression for $\left\|\mathcal{R}(X_j)\right\|$ which allows its efficient evaluation even if $n$ is large.

As a first step, we need to derive an alternative expression for the Riccati residual $\mathcal{R}(X_j)=A^HX_j+X_jA+C^HC-X_jBB^HX_j$. Recall that  $X_j = Z_jY_jZ_j^H$ holds with $Z_j = V_{j+1}\underline{K}_j$ as in \eqref{eq_Z}. Thus, making use of \eqref{rado} we obtain
\[
X_jA =  V_{j+1}\underline{K}_jY_j\underline{K}_j^H V_{j+1}^HA = V_{j+1}\underline{K}_jY_j\underline{H}_j^H V_{j+1}^H.
\]
Moreover, we have $X_jBBH^H_j = V_{j+1}\underline{K}_jY_jS_jY_j\underline{K}_j^HV_{j+1}^H$
with 
\begin{equation}\label{Sj}
S_j = \underline{K}_j^H V_{j+1}^HBB^HV_{j+1}\underline{K}_j = Z_j^H BB^HZ_j.
\end{equation}
 With this, the Riccati residual can be expressed as
\begin{align}
\mathcal{R}(X_j) 
&= V_{j+1}\underline{H}_jY_j \underline{K}_j^HV_{j+1}^H+ V_{j+1}\underline{K}_jY_j \underline{H}_j^HV_{j+1}^H+C^HC- V_{j+1}\underline{K}_jY_j S_jY_j \underline{K}_j^HV_{j+1}^H \nonumber\\
&= V_{j+1}\left( \underline{H}_jY_j \underline{K}_j^H+ \underline{K}_jY_j \underline{H}_j^H+\tilde C \tilde C^H - \underline{K}_jY_j S_jY_j \underline{K}_j^H\right)V_{j+1}^H\label{eq1}
\end{align}
with $C^H = V_{j+1}\tilde C.$ 

Next, rewrite $\Pi_j$ as
\[
\Pi_j = Z_j(W_j^HZ_j)^{-1} W_j^H= V_{j+1}\underline{K}_j(\underline{L}_j^H\underline{K}_j)^{-1} \underline{L}_j^H V_{j+1}^H = V_{j+1}\pi_jV_{j+1}^H
\]
with 
\begin{equation}\label{eq_pismall}
\pi_j \coloneqq \underline{K}_j(\underline{L}_j^H\underline{K}_j)^{-1} \underline{L}_j^H \in \mathbb{C}^{(j+1)p\times (j+1)p}.
\end{equation}
 Note that $\pi_j$ is a projection onto the space spanned by the columns of $\underline{K}_j$ along the orthogonal complement of the space spanned by the columns of $\underline{L}_j.$ This implies $\underline{K}_j^H\pi_j^H = \underline{K}_k^H\underline{L}_j (\underline{K}_k^H\underline{L}_j)^{-1}\underline{K}_j^H = \underline{K}_j^H.$

With this, the projected Riccati residual can be expressed as
\begin{align}
 0 & = \Pi_j \mathcal{R}(X_j)\Pi_j^H \nonumber\\
&= V_{j+1}\pi_j \left( \underline{H}_jY_j \underline{K}_j^H+ \underline{K}_jY_j \underline{H}_j^H+\tilde C \tilde C^H - \underline{K}_jY_j S_jY_j \underline{K}_j^H\right)\pi_j^HV_{j+1}^H.\label{eq2}
\end{align}
As $\Pi_j\mathcal{R}(X_j)\Pi_j^H=0$, \eqref{eq1} and \eqref{eq2} yield the desired alternative expression for the Riccati residual
\begin{align}
\mathcal{R}(X_j)&=\mathcal{R}(X_j)-\Pi\mathcal{R}(X_j)\Pi^H\nonumber\\
&=V_{j+1}\left( \underline{H}_jY_j \underline{K}_j^H+ \underline{K}_jY_j \underline{H}_j^H+\tilde C \tilde C^H - \underline{K}_jY_j S_jY_j \underline{K}_j^H\right)V_{j+1}^H \nonumber\\
&\qquad -  V_{j+1}\pi_j \left( \underline{H}_jY_j \underline{K}_j^H+ \underline{K}_jY_j \underline{H}_j^H+\tilde C \tilde C^H - \underline{K}_jY_j S_jY_j \underline{K}_j^H\right)\pi_j^HV_{j+1}^H \nonumber\\
&= V_{j+1}\left( \tilde \pi_j \underline{H}_jY_j \underline{K}_j^H+ \underline{K}_jY_j \underline{H}_j^H \tilde \pi_j^H +\tilde C\tilde C^H - \pi_j \tilde C \tilde C^H \pi_j^H\right)V_{j+1}^H \label{eq3}
\end{align}
with the projection
 \[
\tilde \pi_j \coloneqq I-\pi_j.
\] 
By construction, $\tilde \pi_j$ is a projection onto the orthogonal complement of the space spanned by the columns of $\underline{L}_j$ along  the space spanned by the columns of $\underline{K}_j.$ That is, 
\begin{equation}\label{eq_tildepi}
\tilde \pi_j = U(W^HU)^{-1}W^H
\end{equation}
with $W\in \mathbb{C}^{(j+1)p\times p}$ such that $\operatorname{range}(W)=\operatorname{range}(\underline{K}_j)^\perp$ and $U\in \mathbb{C}^{(j+1)p\times p}$ such that $\operatorname{range}(U)=\operatorname{range}(\underline{L}_j)^\perp.$

This allows for a formulation of the residual which reduces the storage requirements for evaluation of the residual from an $n \times n$ matrix to a $2p \times 2p$ matrix, independent of the reduced order $(j+1)p$ of \eqref{redric}.
\begin{proposition}\label{theo41}
Let $\Upsilon \coloneqq \underline{K}_jY_j \underline{H}_j^H    + (I-0.5\tilde \pi_j)\tilde C\tilde C^H \in \mathbb{C}^{(j+1)p\times (j+1)p}$ and $T\coloneqq  \Upsilon W (U^HW)^{-1}  \in \mathbb{C}^{(j+1)p\times p}.$ Let $[U~T] = QR,$ with $Q \in \mathbb{C}^{(j+1)p\times 2p}, R \in \mathbb{C}^{2p \times 2p}$
be an economy-size QR decomposition. Then, for any unitarily invariant norm, it holds
\begin{equation}\label{eq_effnorm}
\left\|\mathcal{R}(X_j)\right\|  = \left\| R \left[ \begin{smallmatrix}0 & I_p\\ I_p & 0\end{smallmatrix}\right] R^H\right\|.
\end{equation}
\end{proposition}

\begin{proof}
The constant term in \eqref{eq3} can be written as
\begin{align*}
\tilde C\tilde C^H - \pi_j \tilde C \tilde C^H \pi_j^H 
&= 0.5\left( \tilde C\tilde C^H - \pi_j \tilde C \tilde C^H \pi_j^H - \pi_j \tilde C \tilde C^H  +  \tilde C \tilde C^H \pi_j^H \right) \\
&\qquad +~0.5\left( \tilde C\tilde C^H - \pi_j \tilde C \tilde C^H \pi_j^H + \pi_j \tilde C \tilde C^H  - \tilde C \tilde C^H \pi_j^H \right)\\
&= \tilde \pi_j \tilde C \tilde C^H(I-0.5\tilde \pi_j)^H + (I-0.5\tilde \pi_j)\tilde C\tilde C^H\tilde \pi_j^H
\end{align*}
because of
\begin{align*}
0.5\left( \tilde C\tilde C^H - \pi_j \tilde C \tilde C^H \pi_j^H - \pi_j \tilde C \tilde C^H  +  \tilde C \tilde C^H \pi_j^H \right) 
&=0.5 (I-\pi_j)\tilde C\tilde C^H(I+\pi_j)^H\\
&= 0.5 \tilde \pi_j \tilde C\tilde C^H(I+I-\tilde\pi_j)^H\\
&= \tilde \pi_j \tilde C\tilde C^H(I-0.5\tilde\pi_j)^H,
\end{align*}
and, by a similar argument,
\begin{align*}
0.5\left( \tilde C\tilde C^H - \pi_j \tilde C \tilde C^H \pi_j^H + \pi_j \tilde C \tilde C^H  - \tilde C \tilde C^H \pi_j^H \right)
&= (I-0.5\tilde\pi_j)\tilde C\tilde C^H \tilde \pi_j^H.
\end{align*}
Herewith, \eqref{eq3} can be further manipulated to obtain
\begin{align*}
\mathcal{R}(X_j)&=V_{j+1}\left\{ \tilde \pi_j \left(\underline{H}_jY_j \underline{K}_j^H+  \tilde C \tilde C^H(I-0.5\tilde \pi_j)^H\right) \right.\\
&\qquad\qquad\quad \left.
+~\left(\underline{K}_jY_j \underline{H}_j^H    + (I-0.5\tilde \pi_j)\tilde C\tilde C^H\right)\tilde \pi_j^H   \right\}V_{j+1}^H\\
&=  V_{j+1} \left(\tilde \pi \Upsilon^H +\Upsilon\tilde \pi^H\right) V_{j+1}^H
\end{align*}
with $\Upsilon \coloneqq \underline{K}_jY_j \underline{H}_j^H    + (I-0.5\tilde \pi_j)\tilde C\tilde C^H.$
Hence, for any unitarily invariant norm we have\footnote{The residual formulation \eqref{eq4} is a generalization of the expression in \cite[Proposition 5.3]{Sim16}.}
\begin{equation}\label{eq4}
\left\|\mathcal{R}(X_j)\right\| =\left\|\tilde \pi \Upsilon^H +\Upsilon\tilde \pi^H\right\|.
\end{equation}
 
Making use of \eqref{eq_tildepi} we can further reduce the size of the matrix whose norm has to be computed. It holds
\[
\tilde \pi \Upsilon^H +\Upsilon\tilde \pi^H = U(W^HU)^{-1}W^H\Upsilon^H + \Upsilon W (U^HW)^{-1}U^H = UT^H+TU^H
\]
with $T= \Upsilon W (U^HW)^{-1}.$
As
\[
UT^H +TU^H = [U~T] \left[ \begin{smallmatrix}0 & I_p\\ I_p & 0\end{smallmatrix}\right][U~T]^H
\]
we have with the economy-size QR decomposition
\[
[U~T] = QR,\qquad Q \in \mathbb{C}^{(j+1)p\times 2p}, R \in \mathbb{C}^{2p \times 2p}
\]
for any unitarily invariant norm
\begin{equation}\label{eq_effnorm}
\left\|\mathcal{R}(X_j)\right\|  = \left\| R \left[ \begin{smallmatrix}0 & I_p\\ I_p & 0\end{smallmatrix}\right] R^H\right\|.
\end{equation}
This proves our claim.
\end{proof}

The calculation of the residual norm as described above is summarized in Algorithm \ref{alg2}. There is no incremental update formula how to compute $\left\|\mathcal{R}(X_j)\right\|$ from $\left\|\mathcal{R}(X_{j-1})\right\|,$ but only computations with small-scale matrices of size $(j+1)p\times jp$ and  $(j+1)p\times p$ are involved. 

\begin{algorithm}
\caption{Residual norm calculation}\label{alg2}
\begin{algorithmic}[1]
\Require $\underline{K}_j, \underline{H}_j, \underline{L}_j  \in \mathbb{C}^{(j+1)p \times jp},$ solution $Y_j$ of \eqref{redric} and $\tilde C \in \mathbb{C}^{(j+1)p\times p}$ such that $C^H=V_{j+1}\tilde C.$ 
\Ensure Residual norm $\mathcal{R}(X_j).$
\State  Compute basis $U \in \mathbb{C}^{(j+1)p\times p}$ of $\operatorname{range}(\underline{L}_j)^\perp.$
\State  Compute basis $W \in \mathbb{C}^{(j+1)p\times p}$ of $\operatorname{range}(\underline{K}_j)^\perp.$
\State Compute $\Gamma = W(U^HW)^{-1} \in \mathbb{C}^{(j+1)p\times p}.$
\State Compute $\Psi = \tilde C^H\Gamma \in \mathbb{C}^{(p\times p}.$
\State Compute $T = \underline{K}_jY_j \underline{H}_j^H \Gamma   + (\tilde C-0.5 U\Psi^H) \Psi \in \mathbb{C}^{(j+1)p\times p}.$
\State Compute the economy-size QR decomposition $[U~T] = QR$ with $R \in \mathbb{C}^{2p \times 2p}.$
\State Return $\left\| R \left[ \begin{smallmatrix}0 & I_p\\ I_p & 0\end{smallmatrix}\right] R^H\right\|.$
\end{algorithmic}
\end{algorithm}

\begin{remark}
From 
\[
\mathcal{R}(X_j) =  V_{j+1} \left(\tilde \pi \Upsilon^H +\Upsilon\tilde \pi^H\right) V_{j+1}^H =  V_{j+1} \left(  UT^H+TU^H  \right) V_{j+1}^H 
\]
it follows that the rank of the residual is at most $2p$  as $V_{j+1}$ is of full rank. There are essentially only two scenarios in which the residual rank is smaller. On the one hand the residual can be of rank $p$, which happens for the RADI approximate solution as observed in \cite{BenBKS18}. On the other hand, a rank-$0$ residual is obtained for an exact solution. The rank-$2p$ property has
also been investigated in \cite[Section 5.2]{Wol14} for Lyapunov equations. It was called dilemma of the rational Krylov subspace approach as in general the norm minimizing approximate solution yields a residual with larger rank.
\end{remark}

\section{Truncation of the approximate solution}\label{sec5}
The general projection method generates an approximate solution of the form
\[
X_j= V_{j+1}\underline{K}_jY_j\underline{K}_j^HV_{j+1}^H.
\]
This matrix should be (by construction) Hermitian positive semidefinite with $n-jp$ zero eigenvalues.
But due to rounding errors, it may even be indefinite. We propose to use an eigendecomposition of the Hermitian  matrix $\underline{K}_jY_j\underline{K}_j^H$ in order to truncate non-positive (and possibly some small positive) eigenvalues. This yields a truncated
approximate solution $\widehat X_j,$ which is positive semidefinite  with at least $n-jp$ zero eigenvalues. Thus, the approximate $\widehat{X}_j$ is of lower rank than $X_j.$
This truncated approximate solution can be interpreted as the solution of the Riccati equation projected to a
subspace $\widehat{\mathcal{K}}_j \subset \mathcal{K}_j$ of dimension $r \leq jp$ which is determined by the decomposition of
$X_j.$ Hence, unlike for similar truncated solutions in other projection based approaches, making use of the derivations in Section \ref{sec4}, the residual norm $\|\mathcal{R}(\widehat X_j)\|$ can be evaluated cheaply.

Let 
\begin{equation}\label{brad_here}
A^HV_{j+1}\underline{K}_j = V_{j+1}\underline{H}_j
\end{equation}
 be an orthonormal BRAD as in \eqref{rado} and let $Y_j$ be the solution of \eqref{redric}, so that $X_j= V_{j+1}\underline{K}_jY_j\underline{K}_j^HV_{j+1}^H$ is the current approximate solution to \eqref{ric0}.
Assume that the Hermitian  matrix $\underline{K}_jY_j\underline{K}_j^H \in \mathbb{C}^{(j+1)p \times (j+1)p}$ has been block diagonalized by an unitary matrix such that
\begin{equation}\label{eq7}
\underline{K}_jY_j\underline{K}_j^H
= [\widehat{\underline{Q}}_j~~\check{\underline{Q}}_j] \begin{bmatrix} \widehat Y_j & 0\\ 0 & \check Y_j\end{bmatrix} \begin{bmatrix} \widehat{\underline{Q}}_j^H \\ \check{\underline{Q}}_j^H\end{bmatrix} 
\end{equation}
with $\widehat Y_j = \widehat Y_j^H \in \mathbb{C}^{r \times r},$  $\widehat{\underline{Q}}_j \in  \mathbb{C}^{(j+1)p \times r}$ and $\check{\underline{Q}}_j \in  \mathbb{C}^{(j+1)p \times (j+1)p-r}$ where $\widehat{\underline{Q}}_j^H\widehat{\underline{Q}}_j=I_r, 
\check{\underline{Q}}_j^H\check{\underline{Q}}_j=I_{(j+1)p-r},$ and $\check{\underline{Q}}_j^H\widehat{\underline{Q}}_j=0_{(j+1)p-r\times r}.$
All unwanted ( e.g., all non-positive) eigenvalues of $\underline{K}_jY_j\underline{K}_j^H$ are assumed to be eigenvalues of $\check Y_j.$ Hence, $r \leq jp.$
As the truncated approximate solution define
\[
\widehat X_j = V_{j+1} \widehat{\underline{Q}}_j \widehat Y_j \widehat{\underline{Q}}_j^H V_{j+1}^H.
\]
This can be understood as if $\widehat{X}_j$ is given by its economy-size singular value decomposition.

Next we will show that $\widehat X_j$ is the solution of the Riccati equation \eqref{ric0} projected onto the subspace 
$\widehat{\mathcal{K}}_j = \operatorname{range}(V_{j+1}\widehat{\underline{Q}}_j).$ 
 In order to do, we first derive a BRAD-like expression involving $V_{j+1}\widehat{\underline{Q}}_j.$
Let $T_1 \in \mathbb{C}^{jp\times r}$ be such that $\underline{K}_jT_1 = \widehat{\underline{Q}}_j$ holds. As $\underline{K}_j$ is a matrix of full column rank, its pseudoinverse $\underline{K}_j^+ = (\underline{K}_j^H\underline{K}_j)^{-1}\underline{K}_j^H$ exists. Thus, we have $T_1 = \underline{K}_j^+\widehat{\underline{Q}}_j.$
 Define 
\begin{equation}\label{widehatHj}
\widehat{\underline{H}}_j \coloneqq \underline{H}_jT_1  \in  \mathbb{C}^{(j+1)p \times r}.
\end{equation}
 Then postmultiplying \eqref{brad_here} by $T_1$ yields the BRAD-like relation
\[
A^HV_{j+1}\widehat{\underline{Q}}_j = A^HV_{j+1}\underline{K}_j T_1= V_{j+1}\underline{H}_j  T_1 = V_{j+1}\widehat{\underline{H}}_j.
\]

 Let $\widehat Z_j \coloneqq V_{j+1}\widehat{\underline{Q}}_j \in \mathbb{C}^{n \times r}.$ 
Choose $\widehat{\underline{L}}_j \in \mathbb{C}^{(j+1)p \times r}$ in the same way as $\underline{L}_j$ was chosen in Algorithm \ref{alg1},  that is, e.g., choose $\widehat{\underline{L}}_j = \widehat{\underline{K}}_j,$   $\widehat{\underline{L}}_j = \widehat{\underline{H}}_j,$ or $\widehat{\underline{L}}_j = \widehat{\underline{H}}_j-\widehat{\underline{K}}_j$ in case ${\underline{L}}_j = {\underline{K}}_j,$   ${\underline{L}}_j = {\underline{H}}_j,$ or ${\underline{L}}_j = {\underline{H}}_j-{\underline{K}}_j$ was chosen, resp.. This guarantees by construction that $\operatorname{rank}(\widehat{\underline{L}}_j)=r$, and hence, that
$\widehat{\underline{L}}_j^H\widehat{\underline{Q}}_j$ is nonsingular. 

Let $\widehat W_j \coloneqq V_{j+1}\widehat{\underline{L}}_j \in \mathbb{C}^{n \times r}$ and 
\begin{align*}
\widehat \Pi_j &\coloneqq \widehat Z_j (\widehat W_j^H\widehat Z_j)^{-1} \widehat W_j^H
= V_{j+1} \widehat{\underline{Q}}_j (\widehat{\underline{L}}_j^H\widehat{\underline{Q}}_j)^{-1}\widehat{\underline{L}}_j^HV_{j+1}^H
= V_{j+1} \hat \pi_j V_{j+1}^H
\end{align*}
with $\widehat \pi _j \coloneqq \widehat{\underline{Q}}_j (\widehat{\underline{L}}_j^H\widehat{\underline{Q}}_j)^{-1}\widehat{\underline{L}}_j^H.$ Then $\widehat \Pi_j$
is a (in general oblique) projection onto $\widehat{\mathcal{K}}_j,$ while $\widehat \pi_j$ is a projection on the space spanned by the columns of $\widehat{\underline{Q}}_j.$  Thus, due to \eqref{eq_pismall}, we have $\widehat \pi_j = \widehat \pi_j \pi_j.$

Now we can state and prove the main statement of this section.
\begin{theorem}\label{theo1}
The truncated approximate solution $\widehat X_j = V_{j+1}\widehat{\underline{Q}}_j \widehat Y_j \widehat{\underline{Q}}_j^HV_{j+1}^H$ satisfies the projected equation
\begin{equation}\label{eq5}
\widehat\Pi_j \mathcal{R}(\widehat X_j)\widehat \Pi_j^H =0.
\end{equation}
That is, for $\widehat Y_j \in \mathbb{c}^{r \times r}$ the equation
\begin{equation}\label{eq6}
\widehat \pi_j(\widehat{\underline{H}}_j \hat Y_j \widehat{\underline{Q}}_j^H + \widehat{\underline{Q}}_j \widehat Y_j \widehat{\underline{H}}_j^H + \tilde C \tilde C^H - \widehat{\underline{Q}}_j \widehat Y_j \widehat S_j \widehat Y_j \widehat{\underline{Q}}_j^H )\widehat \pi_j^H=0
\end{equation}
holds, where $\widehat S_j = \widehat{\underline{Q}}_j V_{j+1}^H BB^HV_{j+1}\widehat{\underline{Q}}_j^H .$
\end{theorem}
\begin{proof}
As in \eqref{eq2}, we have for \eqref{eq5}
\begin{align*}
0 &= \widehat\Pi_j \mathcal{R}(\widehat X_j)\widehat \Pi_j^H
= V_{j+1} \widehat \pi_j(\widehat{\underline{H}}_j \hat Y_j \widehat{\underline{Q}}_j^H + \widehat{\underline{Q}}_j \widehat Y_j \widehat{\underline{H}}_j^H + \tilde C \tilde C^H - \widehat{\underline{Q}}_j \widehat Y_j \widehat S_j \widehat Y_j \widehat{\underline{Q}}_j^H )\widehat \pi_j^HV_{j+1}^H.
\end{align*}
This is equivalent to the small scale equation \eqref{eq6}.

 We will now prove that $\widehat{Y}_j$ fulfills \eqref{eq6} which proves the statement of the theorem. In order to do so, we start from the equation
\[
0= \pi_j \left( \underline{H}_jY_j \underline{K}_j^H+ \underline{K}_jY_j \underline{H}_j^H+\tilde C \tilde C^H - \underline{K}_jY_j S_jY_j \underline{K}_j^H\right) \pi_j^H
\]
which is equivalent to \eqref{eq2}.
Pre- and postmultiplication  by $\widehat \pi_j$ and making use of $\widehat \pi_j \pi_j = \widehat \pi_j$ gives
\begin{equation}\label{eq9}
0=\widehat \pi_j \left( \underline{H}_jY_j \underline{K}_j^H+ \underline{K}_jY_j \underline{H}_j^H+\tilde C \tilde C^H - \underline{K}_jY_j S_jY_j \underline{K}_j^H\right)\widehat \pi_j^H.
\end{equation}
We will see that this equation is equivalent to \eqref{eq6}. In order to see this, we consider the different terms one by one.

But first note that due to \eqref{eq7}
\begin{align}
\widehat \pi_j \underline{K}_j Y_j \underline{K}_j^H 
= \begin{bmatrix} \widehat{\underline{Q}}_j & 0\end{bmatrix} \begin{bmatrix}\widehat Y_j &0\\ 0 & \check Y_j\end{bmatrix} 
\begin{bmatrix} \widehat{\underline{Q}}_j^H \\ \check{\underline{Q}}_j^H\end{bmatrix} 
=\widehat{\underline{Q}}_j\widehat Y_j\widehat{\underline{Q}}_j^H 
=\widehat \pi_j \widehat{\underline{Q}}_j\widehat Y_j\widehat{\underline{Q}}_j^H \label{eq10}
\end{align}
holds as $\widehat \pi_j \widehat{\underline{Q}}_j = \widehat{\underline{Q}}_j$ and $\widehat \pi_j \check{\underline{Q}}_j = 0.$

Now consider the first term in \eqref{eq9}. Making use of the transpose of \eqref{eq10} and of \eqref{widehatHj} it holds that
\begin{align*}
\widehat \pi_j  \underline{H}_jY_j \underline{K}_j^H \widehat \pi_j^H 
&= \widehat \pi_j  \underline{H}_j \underline{K}_j^+\underline{K}_j   Y_j \underline{K}_j^H \widehat \pi_j^H
= \widehat \pi_j  \underline{H}_j \underline{K}_j^+\widehat{\underline{Q}}_j   \widehat Y_j \widehat{\underline{Q}}_j^H \widehat \pi_j^H 
= \widehat \pi_j  \widehat{\underline{H}}_j  \widehat Y_j \widehat{\underline{Q}}_j^H \widehat \pi_j^H
\end{align*}
as $\underline{K}_j^+\underline{K}_j =I$ for the pseudoinverse $\underline{K}_j^+.$  Thus, the first term in \eqref{eq9} and \eqref{eq6} are equivalent.

 The second term in \eqref{eq9}  is just the transpose of the first term and is thus equivalent to the second term in \eqref{eq6}.
 The third term in  \eqref{eq9} and \eqref{eq6} are identical.

 For the fourth and last term in \eqref{eq9} it holds with \eqref{Sj} and \eqref{eq10}  that
\begin{align*}
\widehat \pi_j \underline{K}_jY_j S_jY_j \underline{K}_j^H \widehat \pi_j^H 
&=\widehat \pi_j \underline{K}_jY_j \underline{K}_j^HV_{j+1}^HBB^HV_{j+1}\underline{K}_j Y_j \underline{K}_j^H \widehat \pi_j^H\\
&= \widehat \pi_j \widehat{\underline{Q}}_j \widehat Y_j 
\widehat{\underline{Q}}_j^HV_{j+1}^HBB^HV_{j+1}\widehat{\underline{Q}}_j \widehat Y_j \widehat{\underline{Q}}_j^H \widehat \pi_j^H\\
&= \widehat \pi_j \widehat{\underline{Q}}_j \widehat Y_j 
\widehat S_j\widehat Y_j \widehat{\underline{Q}}_j^H \widehat \pi_j^H.
\end{align*}
Thus, the fourth term in \eqref{eq9} and \eqref{eq6} are equivalent. In summary, \eqref{eq9} and \eqref{eq6} are equivalent.
\end{proof}

 In case, 
$\check{Y}_j =0$, \eqref{eq7} reduces to $\underline{K}_jY_j \underline{K}_j^H = \underline{\widehat{Q}}_j\widehat{Y}_j \underline{\widehat{Q}}_j^H.$ Thus, in that case we have $\|\mathcal{R}(X_j)\| = \|\mathcal{R}(\widehat{X}_j)\|$ in any unitarily invariant norm. In any other case, it is not clear whether the norm of $\|\mathcal{R}(\widehat{X}_j)\|$ will decrease or increase compared to $\|\mathcal{R}(X_j)\|.$

Theorem \ref{theo1} allows us to efficiently evaluate the norm of the Riccati residual for the truncated approximate solution 
$\widehat{X}_j$ in the same way as described in the previous section by using $\widehat{\underline{K}}_j  = \widehat{\underline{Q}}_j, \widehat{\underline{H}}_j , \widehat{\underline{L}}_j$ and $\widehat{Y}_j$ instead of the values without $\widehat{~}.$ The solution $\widehat{Y}_j$  of \eqref{eq6}  is given by the decomposition \eqref{eq7}, so no additional small-scale Riccati
equation has to be solved, just a block-diagonalization has to be performed.
The matrices $\widehat{\underline{L}}_j$ and $\widehat{\underline{K}}_j$ consist of $r$ columns, thus we have
$U, W \in \mathbb{C}^{(j+1)p \times p+r}$ in Algorithm \ref{alg2}.
Hence $r < jp$ implies more columns in $U$ and $W$. The rank of the residual matrix increase from $2p$ to $2(p+r)$.  Storing the low-rank factor  $\widehat{Z}_j = V_{j+1}\underline{\widehat{K}}_j$ requires $r$ columns of length $n$, while storing the low-rank factor  ${Z}_j = V_{j+1}\underline{{K}}_j$ requires $jp$ columns of length $n.$ Thus, in case $r < jp$, the low-rank factorization $X_j = \widehat{Z}_j\widehat{Y}_j\widehat{Z}_j^H$ is more efficient storagewise.

In summary, given an approximate $X_j = V_{j+1}\underline{K}_jY_j\underline{K}_j^HV_{j+1}^H$ as discussed in Section \ref{sec3}, one computes \eqref{eq7} and \eqref{widehatHj} and chooses $\underline{\widehat{L}}_j$ appropriately  in order to compute a new approximate $\widehat{X}_j$ of lower rank than $X_j$. Its residual can be determined using the algorithm for the efficient residual norm evaluation (Algorithm \ref{alg2}). The corresponding calculations need to be inserted after line 13  in Algorithm \ref{alg1}. We will term the so modified Algorithm \ref{alg1} the truncation algorithm.

\section{Generalized Riccati Equations}\label{sec6}
For generalized Riccati equations
\begin{equation}\label{ric_gen}
A^HXE+E^HXA+C^HC-E^HXBB^HXE = 0
\end{equation}
with an additional nonsingular matrix $E\in \mathbb{C}^{n \times n}$ it was noted in \cite[Section 4.4]{BenBKS18}, that the equivalent Riccati equation
\[
\mathcal{R}_{gen}(X) = E^{-H}A^HX+XAE^{-1}+E^{-H}C^HCE^{-1}-XBB^HB = 0
\]
has the same structure as \eqref{ric0} where the system matrix $A$ and the initial residual factor $C^H$ are replaced by $AE^{-1}$ and $E^{-H}C^H$, respectively.  In an efficient algorithm inverting $E$ is avoided.
The orthonormal BRAD used in Algorithm \ref{alg1} becomes
\[
E^{-H}A^HV_{j+1}\underline{K}_j = V_{j+1}\underline{H}_h \Leftrightarrow A^HV_{j+1}\underline{K}_j = E^HV_{j+1}\underline{H}_j
\]
with the starting block $E^{-H}C^H.$ The expression \eqref{eq1} for $\mathcal{R}(X_j)$ becomes
\[
\mathcal{R}(X_j) = E^H V_{j+1}\left( \underline{H}_jY_j \underline{K}_j^H+ \underline{K}_jY_j \underline{H}_j^H+\tilde C \tilde C^H - \underline{K}_jY_j S_jY_j \underline{K}_j^H\right)V_{j+1}^HE
\]
with $C^H = E^HV_{j+1}\tilde C$ such that the subsequent derivations in Section \ref{sec3} all hold without any changes.

The expression for the residual becomes 
\begin{equation}\label{EV}
E^HV_{j+1}(UT^H+TU^H)V_{j+1}^HE.
\end{equation}
In order to evaluate this efficiently, let $QR = E^HV_{j+1}\left[ U ~~T\right]$ be an economy-size QR decomposition. Then
\[
\|E^HV_{j+1}(UT^H+TU^H)V_{j+1}^HE\| = \| R \left[\begin{smallmatrix} 0 & I_p\\ I_p&0\end{smallmatrix}\right]R^H\|
\]
is the residual norm of \eqref{ric_gen} (see \eqref{eq_effnorm}).

\section{Numerical Experiments}\label{sec7}
 An extensive comparison of low-rank factored algorithms for solving \eqref{ric0} has been presented in \cite{BenBKS18,BenBKS20}. 
We complement those findings by comparing Algorithm~\ref{alg1} (employing the efficient residual computation as in Algorithm~\ref{alg2} and truncating the approximate solution as discussed in Section \ref{sec5}), the RADI algorithm from \cite{BenBKS18} and the RKSM algorithm from \cite{DruS11,Sim16}  with respect to their convergence performance. Recall that RADI generates approximations $X_k^\text{radi} =  Z_kY_k^\text{radi}Z_k^H$ where the columns of $Z_k$ can be interpreted as a nonorthogonal basis of $\mathcal{K}_k = \mathcal{K}_k(A^H,C^H,\mathcal{S}_k)$ \eqref{eq_ratKryohneC}. The approximations $X_k =  Z_kY_kZ_k^H$ generate via Algorithm \ref{alg1} have the same structure, just the matrix $Y_k$ is obtained in a different way. Moreover, $X_j^\text{radi}$ can be interpreted as the solution of a projection of the large-scale Riccati equation \eqref{ric0} onto the Krylov subspace $\mathcal{K}_k$ employing an oblique projection, see \cite[Section 4.2]{BerF23}. In contrast, RKSM computes an approximate solution $X_k^\text{rksm} =  Q_kY_k^\text{rksm}Q_k^H$ where the columns of $Q_k$ span  an orthogonal basis of $\mathfrak{K}_k = \kappa_k(A^H,C^H,\mathcal{S}_{k-1})$ \eqref{eq_rationalKrylov} and $Y_k^\text{rksm}$ is the solution of  \eqref{smallVk} resulting from the Galerkin projection $Q_k^H \mathcal{R}(X_k) Q_k = 0$ of the Riccati residual.

The behavior of Algorithm~\ref{alg1} for the choices $\underline{L}_j = \underline{K}_j$ (which yields a Galerkin projection) and $\underline{L}_j = \underline{H}_j$ or $\underline{L}_j = \underline{H}_j - \underline{K}_j$ (which yield a Petrov-Galerkin projection),  and different shift strategies is investigated. That is, in order to see the influence of different shift strategies on Algorithm \ref{alg1}, a set of shifts is precalculated. Then all three algorithms are run with this same set of shifts so that all algorithms perform the same number of iteration steps with the same shifts. Their convergence behavior is compared. Our purpose is not to propose Algorithm \ref{alg1} in its current version as a valid competitor of RADI or RKSM. In theory, the Petrov-Galerkin approach gives more degrees of freedom than the Galerkin approach. With an optimal choice of $L_k$ and suitable shifts, the Petrov-Galerkin version of Algorithm \ref{alg1} should converge at least as fast or faster than the Galerkin version. But so far, we do not know how to best choose $L_k$ and the set of shifts.

For the RADI algorithm we use the function \texttt{mess\_lrradi} from the MATLAB toolbox M.E.S.S.-2.2 \cite{SaaKB22}. The implementation of the RKSM method used is based on the function  \texttt{RKSMa\_care}  from \cite{BenBKS20}. The function \texttt{RKSMa\_care} had to be modified in order to handle precomputed shifts correctly. The modified code as well as all experimental code used to generate the results presented can be found at \cite{zenodo}.

We consider three different shift strategies. For two of these, \texttt{mess\_lrradi} is used.} The shifts are precomputed employing either the default  option \texttt{opts.shifts.method = 'gen-ham-opti';}   which implies that the shift strategy residual Hamiltonian shifts \cite[Section 4.5.2]{BenBKS18} is used or the option \texttt{opts.shifts.method = 'heur';} which yields an estimation of suboptimal ADI shift parameters  as suggested by Penzl \cite{lyapack, Kue16}. The third set is computed using the default option \texttt{conv} of the function \texttt{RKSMa\_care} which chooses the shifts from convex hull of Ritz values \cite{DruS11}.

Algorithms \ref{alg1} and \ref{alg2} have been  implemented in order to handle generalized Riccati equations as discussed in Section \ref{sec6}. 
The entire orthonormal BRAD \eqref{rado} associated to the precomputed set of shifts is generated prior to the start of the iteration in Algorithm~\ref{alg1} using the function \texttt{rat\_krylov} from the Rational Krylov Toolbox \cite{BerEG14} in version 2.9. During the iteration, only the relevant parts of the BRAD are used. The \texttt{rat\_krylov} function supports generalized Riccati equations with an additional system matrix $E$, block vectors
$C^H \in \mathbb{C}^{n\times p}$ and realification in case complex shifts are used in conjugate-complex pairs. All occurring small-scale Riccati equations are solved with MATLABs \texttt{icare}. For the vast majority of the small-scale Riccati equations to be solved, MATLAB's \texttt{icare} reported back that the unique solution generated is accurate (\texttt{info.Report == 0}).  Nonetheless, some of the computed solutions $Y_j$ were indefinite. In our implementation $\check{Y}_j$ in \eqref{eq7} was chosen to contain all eigenvalues less than  $10^{-12}\rho(\underline{K}_jY_j\underline{K}_j^H)$ where $\rho(\cdot)$ denotes the spectral radius of  the matrix. In order to do so, the matrix $\underline{K}_jY_j\underline{K}_j^H$ is diagonalized using MATLAB's \texttt{eig}. The eigenvalues (and corresponding eigenvectors) are reordered such that the eigenvalues appear in descending order.

Algorithm \ref{alg1} has been implemented in two versions: an efficient version that takes into account the comments from Section \ref{subsec_eff} as well as a version that directly implements Algorithm \ref{alg1}. The efficient version of the Algorithm \ref{alg1} does require $\underline{K}_j,$ $\underline{H}_j$ and $\underline{L}_j$ to be block upper Hessenberg matrices.  The function \texttt{rat\_krylov} may return $\underline{K}_j$ and $\underline{H}_j$ with additional entries in case deflation was performed. In those cases, we make use of the second version of the implementation of Algorithm \ref{alg1}. In a really efficient implementation, the calculation of the BRAD would have to be linked to the successive solution of the Riccati equation, so that the BRAD is not calculated in advance, but step by step immediately followed by the solution of the corresponding small Riccati equation. Aspects such as deflation and re-orthogonalization would then be adapted to the problem at hand. But for our purpose, the approach taken here is sufficient.

In our test set up, RADI will, in general, be (much) faster than RKSM and Algorithm \ref{alg1} as each algorithm is run for the same number of iteration steps with the same set of shifts. The most timing consuming part of the RADI algorithm is the solution of linear systems with multiple right-hand sides.  The rest of the computational effort is neglectable. Assuming that, as proposed in  \cite[Section 4.2]{BenBKS18}, the Sherman-Morrison-Woodbury formula is used to reformulate the dense linear systems in the RADI algorithm as sparse linear systems of the form $(A^H-\sigma I)\hat Z = R,$ linear systems with $p+m$ right-hand sides have to be solved. RKSM and Algorithm~\ref{alg1} need to solve the same linear systems of equations with just $p$ right hand sides in order to advance the required basis of the block rational Krylov subspace. However, this slight advantage is more than offset by the necessary orthonormalization of the basis in RKSM and Algorithm~\ref{alg1}.
 Moreover, unlike in the RADI algorithm, in RKSM and Algorithm~\ref{alg1} the solution of the small-scale Riccati equation \eqref{redric} of growing dimension has to be calculated, which further increases the computational time required. 
Nevertheless, for the first example we present time measurements to illustrate the difference between RKSM, Algorithm \ref{alg1} and RADI and especially between RKSM and Algorithm \ref{alg1}.

All experiments are performed in MATLAB R2024a on an Intel(R) Core(TM) i7-8565U CPU @ 1.80GHz
1.99 GHz with 16GB RAM.

\subsection{Example 1}
The first example considered is the well-known steel profile cooling model
from the Oberwolfach Model Reduction Benchmark Collection \cite{morwiki_steel,BenS05b}. This example (often termed RAIL)  comes in different problem sizes $n$, but fix $m =7$ and $p = 6$. We used the one with $n = 79,841.$ The system matrices $E$ and $A$ are symmetric positive and negative definite, resp.. 

\begin{figure}[!htb]
\begin{center}
\includegraphics[scale=0.45]{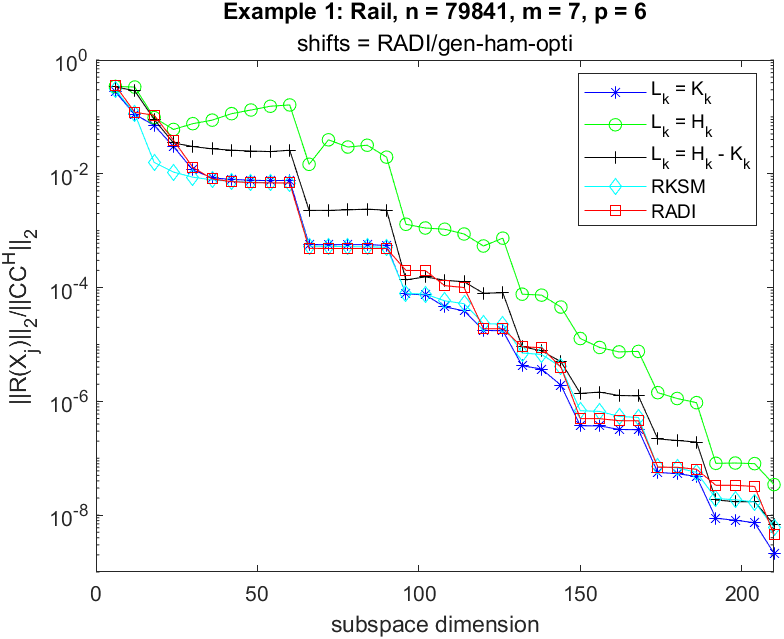} \includegraphics[scale=0.45]{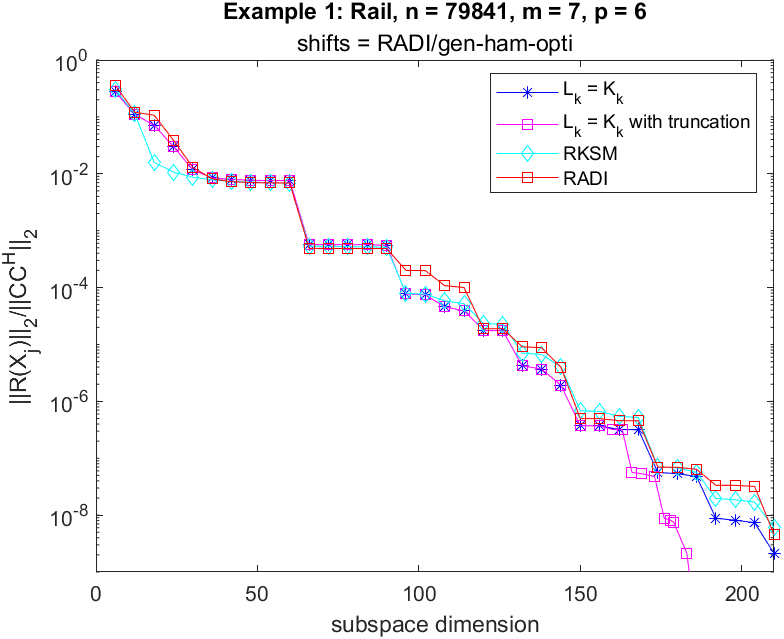}

\includegraphics[scale=0.45]{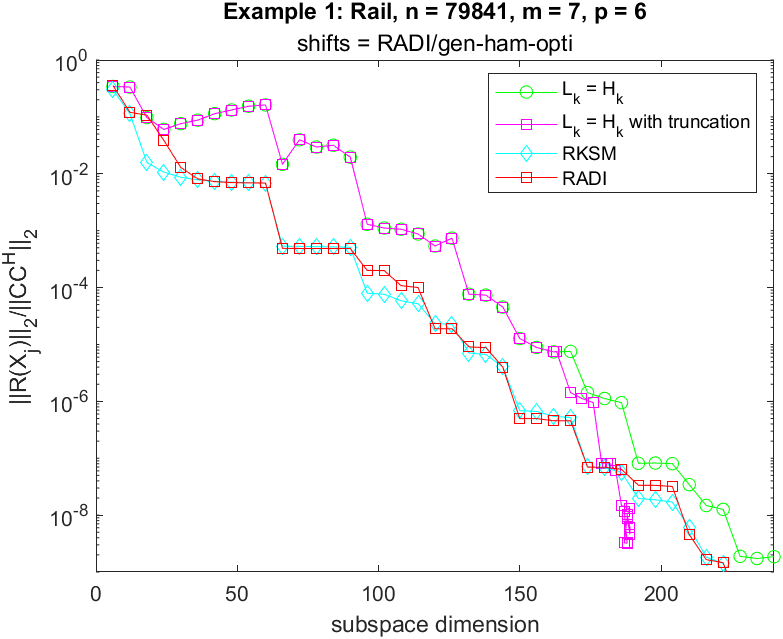}  \includegraphics[scale=0.45]{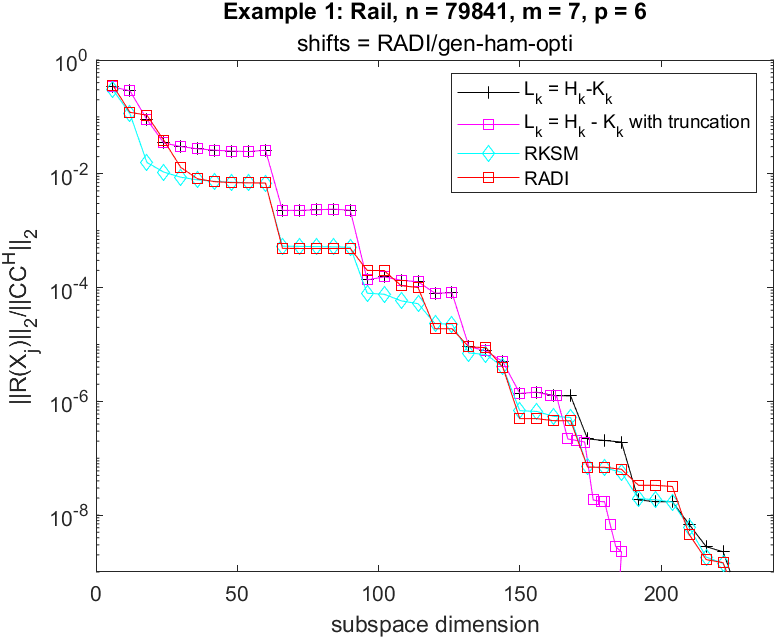}
\end{center}
\caption{Example 1: Relative residual norms for shifts generated by \texttt{mess\_lrradi} with option \texttt{'gen-ham-opti'}. \label{fig1}}
\end{figure}

As a first test we precomputed  one set of shifts employing \texttt{mess\_lrradi} with the option \texttt{'gen-ham-opti'}.   For all other options used, we refer to \cite{zenodo}. This resulted in $50$  real shifts.  Then we ran all algorithms with this set of shifts. All plots in Figure \ref{fig1} display the residual norm normalized by the norm of $CC^H$ versus the dimension of the subspace projected onto for all methods considered. The plot in the upper left-hand side shows that Algorithm \ref{alg1} with the  choice $\underline{L}_j = \underline{H}_j$ and $\underline{L}_j=\underline{H}_j-\underline{K}_j$ converges a bit slower than RADI and RKSM. For the choice $\underline{L}_j=\underline{K}_j$ Algorithm \ref{alg1} performs slightly better than RADI and (in the end) than RKSM.
The other plots in Figure \ref{fig1} show the effect of the truncation on the convergence. To make the effect clearly visible, in each of those two plots three curves are given  that were already shown in the plot at the left in Figure \ref{fig1}, the convergence of RADI, RKSM  and that of Algorithm 1 for one of the three choices of $\underline{L}_j.$ 
In addition, the convergence of the truncation algorithm is given for the selected $L_k$. The  truncated approximate solution $\widehat{X}_j$ is computed from $X_j$ as explained in Section \ref{sec5}. For illustration purposes this is done in in each iteration step. Clearly, in  the beginning no truncation is taking place, but after $26$ iteration steps (that is, BRAD subspace dimension $162$) truncation shows an effect. As can be seen in Figure \ref{fig1} and Table \ref{tab2}, when employing truncation of the computed approximate solution ${X}_j,$ the subspace dimensions decrease after a while such that the low-rank factorization of $\widehat{X}_j$ requires less storage than that of $X_j.$ Recall that when applying Algorithm \ref{alg1}, the dimension of the subspace on which the Riccati equation is projected  increases by $p$ in each iteration step, that is, in step $j$, it is $6j.$ Moreover, in each iteration step, the rank of the Riccati residual $\mathcal{R}(X_j)$ is given by $2p = 12.$ See Table \ref{tab2} for detailed information on the dimension of the subspace used and the rank of the Riccati residual. In the end, instead of an $79,841 \times 300$ matrix to store $V_{j+1}\underline{K}_j$ only one of size $79,841 \times 187$ / $79,841 \times 188$ / $79,841 \times 188$ required, depending on the choice of $\underline{L}_j.$

\begin{table}[htb!]\footnotesize
\begin{center} 
\begin{tabular}{|c|c|c|c|c|c|c|}\hline
$j$ &\multicolumn{2}{c|}{$\underline{L}_j = \underline{K}_j$} & \multicolumn{2}{c|}{$\underline{L}_j = \underline{H}_j$} & \multicolumn{2}{c|}{$\underline{L}_j = \underline{H}_j-\underline{K}_j$}\\ \hline
&  & subspace &   & subspace & & subspace \\ 
& rank$(\mathcal{R}(X_j))$ & dimension     &rank$(\mathcal{R}(X_j))$ & dimension & rank$(\mathcal{R}(X_j))$ & dimension\\ \hline
 $1-26$  & $12$ & $6j$ &  $12$  & $6j$ &   $12$ & $6j$\\
 $27$ &  $16$ & $160$   & $12$   & $162$ & $14$  & $161$\\
$28$ &  $22$  & $163$  & $22$    & $163$ & $22$   & $163$\\
$29$ &  $28$  & $166$  & $24$    & $168$ & $26$  & $167$ \\
$30$ &  $34$  & $169$  & $28$    & $172$ & $32$  & $170$\\
$31$ &  $38$  & $173$  & $32$    & $176$ & $38$  & $173$ \\
$32$ &  $44$  & $176$  &  $38$   & $179$ & $44$  & $176$\\
$33$ &  $52$  & $178$  &  $46$   & $181$ & $50$  & $179$\\
$34$ &  $62$  & $179$  & $56$    & $182$ & $60$  & $180$\\  
$35$ &  $66$  & $183$  &  $64$   & $184$ & $68$  & $182$\\
$36$ &  $76$  & $184$  &  $72$   & $186$ &  $76$ & $184$\\
$37$ &  $88$  & $184$  &  $78$   & $189$ &  $84$ & $186$\\
$38$ &  $98$  & $185$  &  $92$   & $188$ & $98$  & $185$\\
$39$ &  $110$ & $185$  & $102$  & $189$ & $106$ & $187$\\
$40$ &  $120$ & $186$  & $114$  & $189$ & $118$ & $187$\\
$41$ &  $132$ & $186$  & $128$  & $188$ & $130$ & $187$\\
$42$ &  $144$ & $186$  & $140$  & $188$ & $142$ & $187$\\
$43$ &  $154$ & $187$  & $154$  & $187$ & $154$ & $187$\\
$44$ &  $166$ & $187$  & $166$  & $187$ & $166$ & $187$\\
$45$ &  $178$ & $187$  & $176$  & $188$ & $176$ & $188$ \\
$46$ &  $190$ & $187$  & $188$  & $188$ & $188$ & $188$\\
$47$ &  $202$ & $187$  & $200$  & $188$ & $200$ & $188$\\
$48$ &  $214$ & $187$  & $212$  & $188$ & $212$ & $188$\\
$49$ &  $226$ & $187$  & $224$  & $188$ & $214$ & $188$\\
$50$ &  $238$ & $187$  & $236$  & $188$ & $238$ & $188$\\ \hline
\end{tabular}
\end{center}
\caption{Example 1: Rank$(\mathcal{R}(X_j))$ and subspace dimension for the set of \texttt{'gen-ham-opti'} shifts as in Figure \ref{fig1}.\label{tab2}}
\end{table}

\begin{figure}[htb!]
\begin{center}
\includegraphics[scale=0.45]{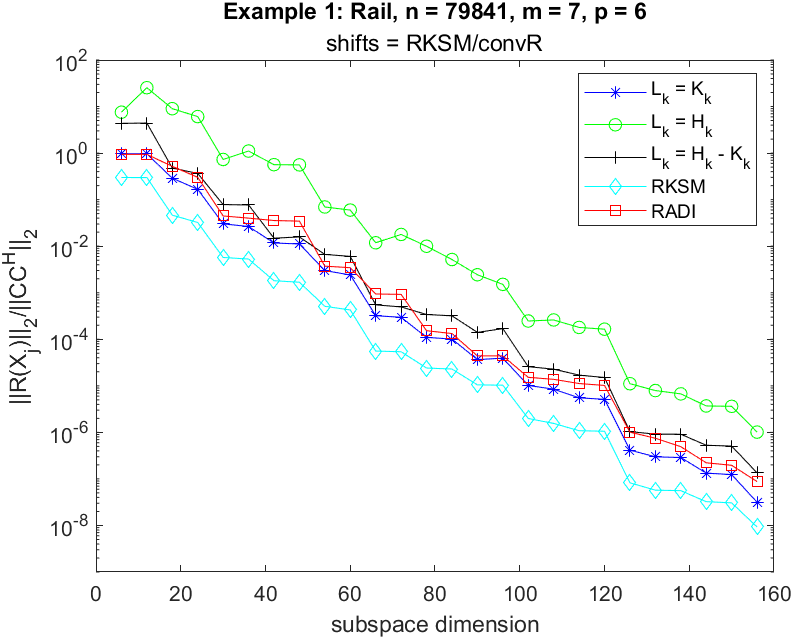}  \includegraphics[scale=0.45]{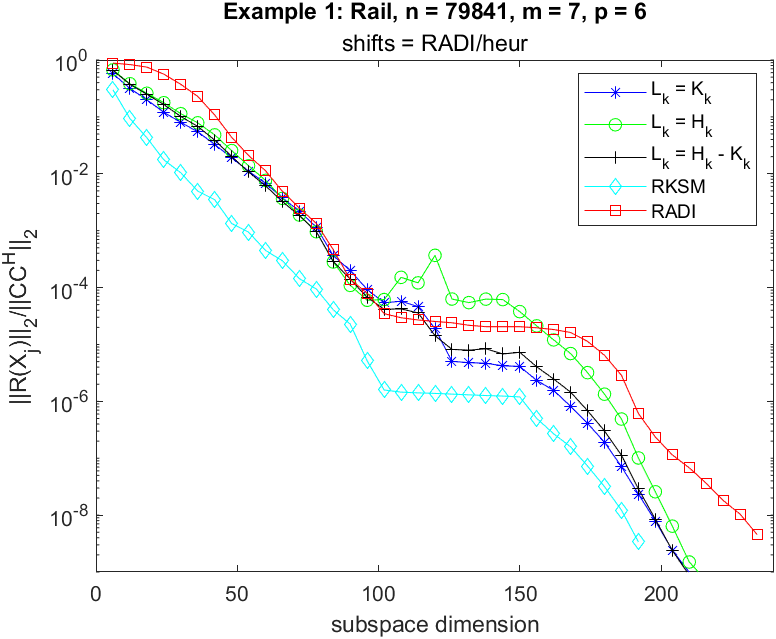}
\end{center}
\caption{Example 1: Relative residual norms for different shifts. \label{fig2}}
\end{figure}

\begin{figure}[htb!]
\begin{center}
\includegraphics[scale=0.45]{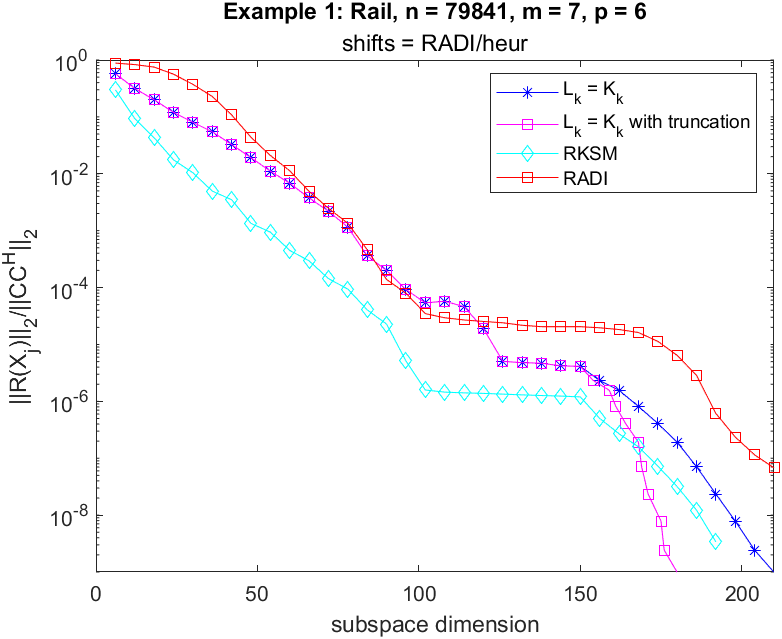}  \includegraphics[scale=0.45]{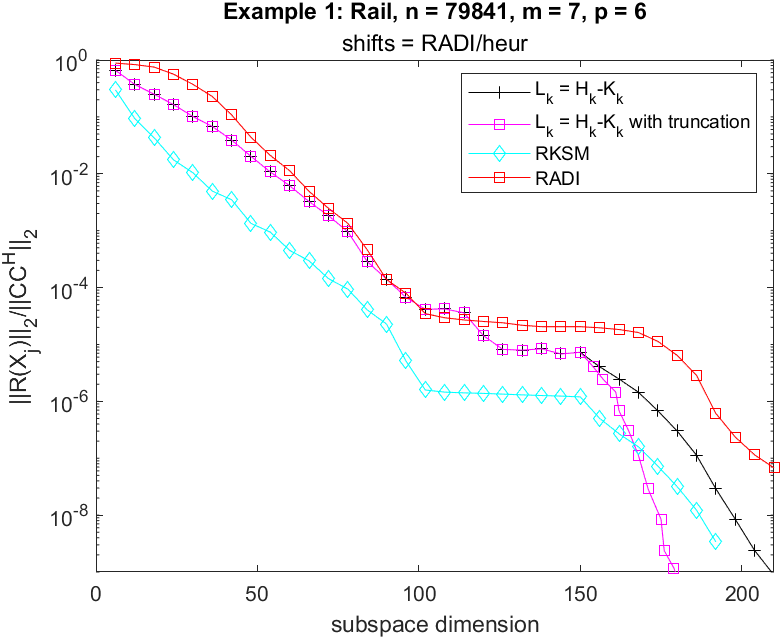}
\end{center}
\caption{Example 1: Relative residual norms for shifts generated by \texttt{mess\_lrradi} and the option \texttt{ 'heur'}. \label{fig3}}
\end{figure}

Figure \ref{fig2} displays the same information as the left upper plot in Figure \ref{fig1}, but with different precomputed shifts. The left plot shows the convergence related to the shifts generated with \texttt{my\_RKSMa\_care} and the  option \texttt{'convR'}, the right one to the convergence related to the shifts generated with \texttt{mess\_lrradi} and the option \texttt{ 'heur'}. 
All precomputed shifts are real. For both set of shifts, RKSM performs better than the other algorithms.
Algorithm \ref{alg1} with the  choice $\underline{L}_j = \underline{K}_j$  converges slightly faster than RADI, but a bit slower than RKSM. For the set of RKSM shifts, the other versions of Algorithm 1 perform worse than the other algorithms, but for the set of \texttt{'heur'} shifts, the convergence for the choice  $L_k=H_k-K_k$ is better than RADI and almost alike the choice $L_k=K_k$,
while for the choice $\underline{L}_j=\underline{H}_j$ the overall performance is better than that of RADI although this does not apply to some iteration steps. Truncation has almost no effect for the set of RKSM shifts, while for the set of \texttt{'heur'} shifts the effect can be nicely seen from BRAD subspace dimension $162$ onwards, see Figure \ref{fig3}.

Finally, we give some timings for the different algorithms and sets of shifts, see Table \ref{tab1}. Algorithm 1 is used with the efficient implementation discussed in Section \ref{subsec_eff} without any truncation. The number of shifts is chosen such that the algorithms terminate with a comparable final accuracy. Due to the set up of the experiments, RADI is the fastest algorithm as both other algorithms need to perform some orthogonalization to set up the required basis. As the timings show, Algorithm \ref{alg1} requires typically less time than RKSM.  

\begin{table}[htb!]\footnotesize
\begin{center} 
\begin{tabular}{|c|c|c|c|}\hline
Shift Choice & Method & Timings & $\sharp$ of shifts\\ \hline
& RKSM &20.1088&\\ \cline{2-3}
& RADI & 12.1011&\\ \cline{2-3}
RKSM& Alg.1,$ L_k=K_k$ &15.8956& 26 real \\\cline{2-3}
& Alg.1, $L_k=H_k$& 15.8769& \\\cline{2-3}
&Alg.1, $L_k = H_k-K_k$& 15.8989&\\ \hline
& RKSM &31.6462 &\\ \cline{2-3}
& RADI & 12.4941 &\\ \cline{2-3}
\texttt{'gen-ham-opti'}&Alg.1, $L_k=K_k$ &26.7037 & 32 real \\\cline{2-3}
& Alg.1, $L_k=H_k$&26.8429 & \\\cline{2-3}
&Alg.1, $L_k = H_k-K_k$&27.6092& \\ \hline
& RKSM &49.8927&\\ \cline{2-3}
& RADI &18.8875 &\\ \cline{2-3}
\texttt{'heur'}&Alg.1, $L_k=K_k$ & 26.5591& 32 real\\\cline{2-3}
&Alg.1, $L_k=H_k$&26.2650&\\\cline{2-3}
&Alg.1, $L_k = H_k-K_k$&26.4021& \\ \hline
\end{tabular}
\end{center}
\caption{Example 1: Timings for different sets of shifts.\label{tab1}}
\end{table}

\subsection{Example 2}
The second example considered is the convection-diffusion benchmark example from MORwiki - Model Order Reduction Wiki \cite{morwiki,lyapack}. The examples are constructed with 

\begin{figure}[!htb]
\begin{center}
\includegraphics[scale=0.45]{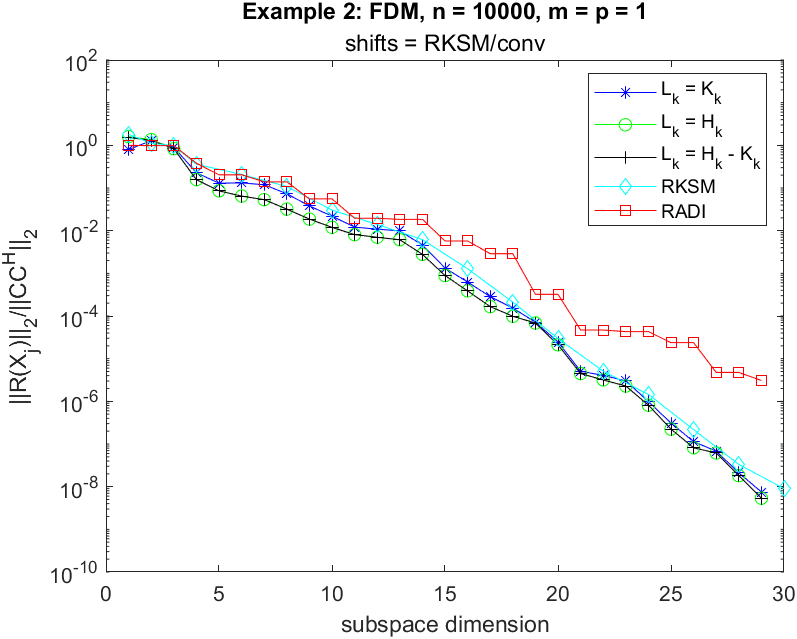} \includegraphics[scale=0.45]{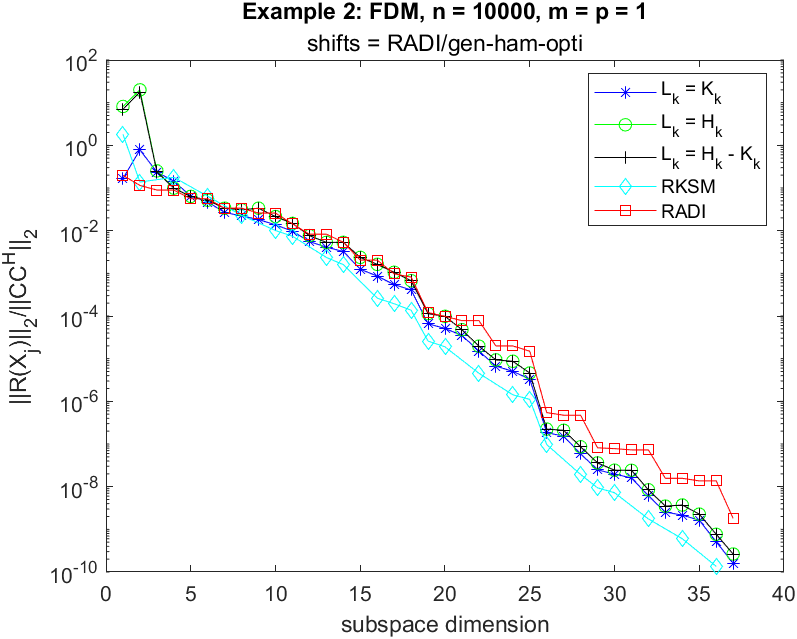}  
\end{center}
\caption{Example 2: Relative residual norms for different set of shifts.\label{fig4}}
\end{figure}

\begin{center}
\begin{minipage}{12cm}
\texttt{A = fdm\_2d\_matrix(100,'10*x','100*y','0');}\\
\texttt{B = fdm\_2d\_vector(100,'.1<x<=.3');}\\
\texttt{ C = fdm\_2d\_vector(100,'.7<x<=.9')'; }\\
\texttt{E = speye(size(A));} 
\end{minipage}
\end{center}

\noindent
resulting in a SISO system of order $n = 10,000.$  The plots in Figure \ref{fig4}  display essentially  the same information as the corresponding plots in Figure \ref{fig2}. For this example, for the set of shifts generated with \texttt{my\_RKSMa\_care} and the  option \texttt{'conv'} all versions of Algorithm \ref{alg1} perform at least as good as RKSM, while for the set of shifts generated with \texttt{mess\_lrradi} and the  option \texttt{'gen-ham-opti'} the performance of all versions of Algorithm \ref{alg1} lie between the one of RKSM and RADI. Truncation has essentially no effect for the RKSM shifts, while for the \texttt{'gen-ham-opti'} shifts, starting from the BRAD subspace dimension 25, truncation has an effect.

\begin{figure}[htb!]
\begin{center}
\includegraphics[scale=0.45]{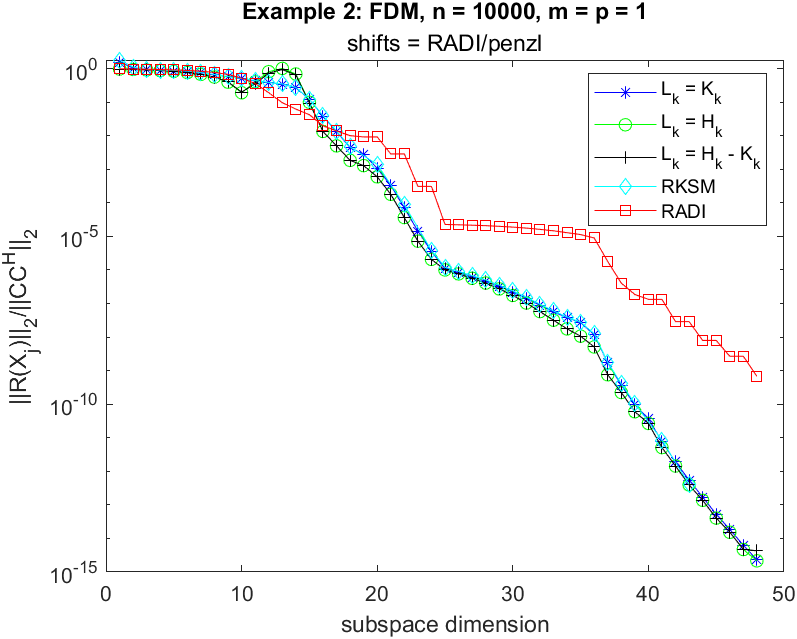}  \includegraphics[scale=0.45]{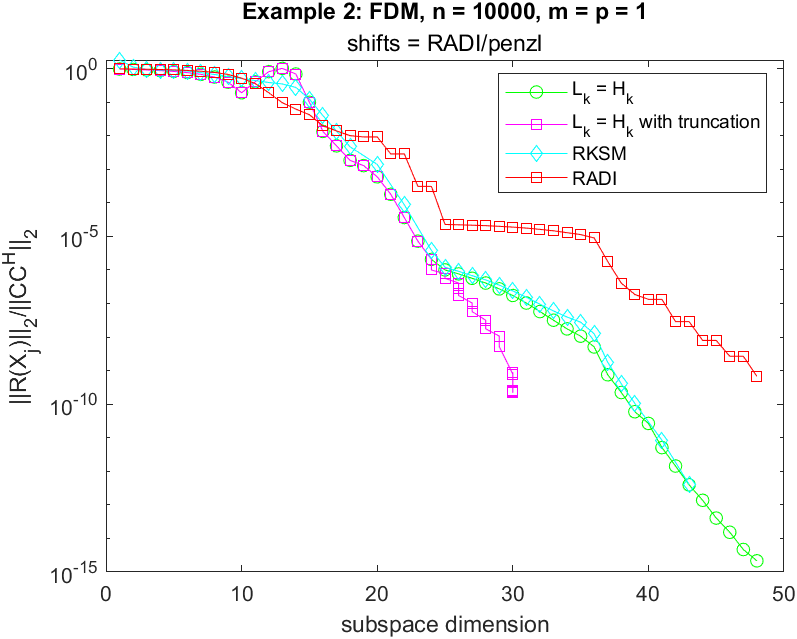}
\end{center}
\caption{Example 2: Relative residual norms for  \texttt{'heur'} shifts generated with \texttt{mess\_lrradi}. \label{fig5}}
\end{figure}

For the set of shifts generated with \texttt{mess\_lrradi} and the  option \texttt{'heur'} RADI performs worse than the other algorithms which perform alike, see Figure \ref{fig5}. Truncating the solution does decrease the final subspace dimension, but also decreases the accuracy of the computed solution. This is depicted here for the choice $L_k = H_K.$

\subsection{Example 3}
\begin{figure}[!htb]
\begin{center}
\includegraphics[scale=0.45]{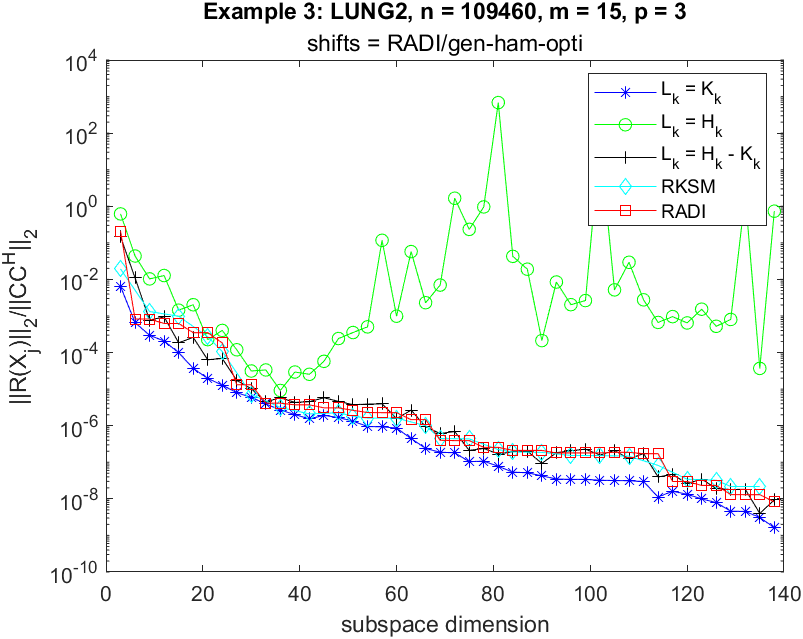} 
\end{center}
\caption{Example 3: Relative residual norms for shifts generated using the option \texttt{'gen-ham-opti'} for \texttt{mess\_lrradi}.\label{fig6}}
\end{figure}

As a third example we consider the nonsymmetric matrix lung2 available from The SuiteSparse Matrix Collection\footnote{\url{https://sparse.tamu.edu}}  (formerly known as the University of Florida Sparse Matrix Collection) via the matrix ID 894 \cite{DavH11}, modeling processes in the human lung. 
We employ this example with the negated system matrix $-A \in \mathbb{R}^{109460 \times 109460}, E = I$ and randomly chosen $C^H \in \mathbb{R}^{109460 \times 3}, B \in \mathbb{R}^{109460\times 15}$ (using \texttt{randn}). Here we report only on the numerical experiments involving the precomputed shifts using the option \texttt{'gen-ham-opti'} for \texttt{mess\_lrradi}. 
When using the shifts precomputed with the \texttt{'heur'} option, MATLAB's \texttt{icare} had trouble solving the small-scale Riccati equations. RKSM also had problems converging.
Figure \ref{fig6} displays the same information as Figure \ref{fig1}. As can be seen Algorithm \ref{alg1} with the choice $\underline{L}_k = \underline{K}_j$ outperforms RKSM and RADI slightly, while for the choice $\underline{L}_k = \underline{H}_j -\underline{K}_j$  RKSM and RADI performs slightly better.  Algorithm \ref{alg1} with the choice $\underline{L}_k = \underline{H}_j$ behaves quite different than the other two choices for $\underline{L}_j.$  For this choice only erratic and poor convergence can be observed.

\subsection{Summary of Findings}
The convergence of the general projection method typically behaves comparable to that of the RADI  and the RKSM algorithm. The general projection method with orthogonal projection ($\underline{L}_j = \underline{K}_j$) shows the best convergence behavior among the  all variants tested, usually slightly better than the RADI algorithm. This may be due to the fact, that all shift choices tested are stemming from Galerkin projection type methods. A choice adapted to the Petrov-Galerkin scenario could possibly provide a remedy and better performance for the Petrov-Galerkin projection approach. The general projection method is computationally more demanding than the RADI algorithm, but somewhat cheaper than RKSM. Truncation of the approximate solution $X_j$  proves to be an efficient way to further reduce the rank of the approximate solution while improving the accuracy of the approximate solution achieved for the related subspace dimension. The choice of criterion for determining $\check Y_j$ \eqref{eq7} determines the achievable accuracy of the truncated approximate solution.

\section{Concluding Remarks}
So far, the orthogonal projection onto the block rational Krylov subspace $\mathfrak{K}_j$ \eqref{eq_rationalKrylov} or onto the extended block Krylov subspace $\kappa_j(A^H,C^H) + \kappa_j(A^{-H},A^{-H}C^H)$ for the standard block Krylov subspace $\kappa_j(A^H,C^H)$ \eqref{eq_standardKrylov}
 has been considered in the literature in connection with projection methods for solving large Riccati equations. Here, for the first time, we have explicitly considered the projection of the Riccati equation onto the block rational Krylov subspace $\mathcal{K}_j$ \eqref{eq_ratKryohneC}. The projections need not be orthogonal. Like the resulting projected small-scale Riccati equation,  the projections are determined by the matrices in the BRAD corresponding to $\mathfrak{K}_j$. An efficient way to evaluate the norm of the residual has been derived. Instead of the norm of an $n \times n$ matrix, the norm of a readily available $2p \times 2p$ matrix has to be computed. This implies that the rank of the residual matrix is $2p.$ The idea of truncating the resulting approximate solution has been proposed. By employing this idea, the rank of the approximate solution can be effectively reduced further while increasing the accuracy obtained for the corresponding subspace dimension.
It has been proven that the  truncated approximate solution can be interpreted as the solution of the Riccati residual projected to a subspace of the Krylov subspace $\mathcal{K}_j.$ This gives a way to efficiently evaluate the norm of the resulting residual. Numerical experiments demonstrate that the convergence of the proposed projection methods generally follows the same pattern as the convergence of the RADI and the RKSM algorithm. Among all the evaluated versions, the general projection method with orthogonal projection exhibits the best convergence behavior. The Petrov-Galerkin scenario allows for more freedom in the choice of the search space. This could possibly provide a remedy and better performance for the Petrov-Galerkin projection approach when a suitable strategy for choosing $\underline{L}_j$ and the set of shifts has been found. However, the RADI approach requires less computational time than the projection method.  

\section*{Acknowledgements}
Part of this work was done while the second author visited the Oden Institute at the University of Texas at Austin in October 2023 and the Department of Mathematics and the Division of Computational Modeling and Data Analytics (CMDA) in the College of Science at Virginia Tech in Blacksburg in November 2023.


\end{document}